\documentclass[preprint,10pt,5p,twocolumn,authoryear]{elsarticle}
\usepackage{amssymb}
\usepackage{amsmath}
\usepackage[colorlinks]{hyperref}
\usepackage{comment}
\AtBeginDocument{%
  \hypersetup{
    citecolor=red,
    linkcolor=blue,   
    urlcolor=black}}
\begin{document}
\begin{frontmatter}
\title{Convex Formulation of the Maritime Fleet Size and Mix Problem Considering Battery Electric Ships}
\author{Antti Ritari\corref{cor1}}
\cortext[cor1]{Corresponding author}
\ead{antti.ritari@aalto.fi}
\author{Jani Romanoff} 
\author{Kari Tammi} 
\affiliation{organization={Department of Energy and Mechanical Engineering, Aalto University},
            city={Espoo},
            country={Finland}}
\begin{abstract}
This paper focuses on the problem of determining a minimum-cost fleet of battery electric ships for a given liner shipping operation. 
The problem is strongly nonlinear and includes integer-valued decision variables, which make it intractable for most real-world instances. The conventional approach in the literature is to formulate a linear approximation by restricting available ship types to a small number of predetermined alternatives. Contrary to the conventional linearization approach, this paper models the nonlinearities directly. We show that the problem exhibits a hidden convex structure uncovered by changes of variables.
Computational experiments show that our convex formulation achieves 15.1\% lower fleet cost on average and reduces the solve time by at least 10x compared to the linear formulation. The solve time advantage is attributed to the elimination of binary variables associated with ship type selection.
Our implementation is available at \url{https://users.aalto.fi/~aritari/convex_ESFSMP.html}.
\end{abstract}

\end{frontmatter}
\tableofcontents
\clearpage
\section{Introduction} \label{sec:introduction}
Direct electrification is emerging as the most economically and technically viable zero-emission pathway for short-sea liner shipping \citep{KPP22, KSG+25}. Direct electrification refers to storing electrical energy onboard the ship using battery energy storage. The well-to-wake efficiency of the direct electrification pathway is at least four times greater compared to the electrofuel pathway \citep{UBD+21}. However, the current commercially available lithium-ion cells exhibit at least two orders of magnitude lower volumetric energy density than liquid hydrocarbon fuels \citep{Pan24}. The low energy density of batteries, coupled with increased charging time compared to liquid fuel bunkering, requires restructuring of most short-sea shipping services. 

Strategic planning a battery-electric liner shipping service presents a complex optimization problem due to the strong coupling of decisions concerning cargo routing, sailing schedules, onshore charger powers, ship sizing, and battery capacity. Accounting for all the details is beyond the capacity of human planners relying on elementary spreadsheet calculations. Therefore, advanced decision support tools building on mathematical optimization modeling can be of great use. The ro-ro, ro-pax, and container feeder segments constitute over 1,350 ships (over 5,000 GT) in total in the European Union waters alone \citep{Ep25}. Clearly, intelligent planning of electrifying these ships will produce large savings and accelerate emission reduction goals. 

This paper focuses on a particular strategic planning problem in the ro-ro, ro-pax, and container feeder liner shipping segments: determining a plan to purchase a fleet of battery-electric ships for a given operation. Our starting point is a group of predetermined shipping routes operated by a liner shipping company. The decision problem is to determine the size, type, and number of ships to purchase and deploy on each route. This fleet design decision is coupled with tactical decisions concerning the number of round-trip voyages each ship completes during the planning horizon. The objective is to minimize the total cost of fleet purchase and operation, subject to a minimum transportation service level constraint. We label this problem the \textit{electric ship fleet size and mix problem} (ESFSMP).

The literature on ESFSMP and its extension to a multi-period setting \citep{HLN+22, HFN+24, CFJ+25} formulate mixed-integer linear programming-based models. These works model the nonlinear relationship between speed and energy consumption by linear interpolation between discrete speed and energy levels. They also introduce binary decision variables for selecting ship types from a predetermined set of alternatives. The linear modeling approach comes with three downsides:
\begin{enumerate}
    \item Energy consumption and transit time are overestimated for all speeds that deviate from the discrete levels, as discussed in detail in \cite{AFH15}.
    \item Constructing the set of ship types requires case-specific expert knowledge. Solution quality suffers if the optimal ship type is excluded from the set.
    \item Increasing the set of ship types increases not only the number of ship type binary variables but also the number of all other discrete and continuous variables that are indexed by ship type. The combinatorial search space grows exponentially, risking intractability.
\end{enumerate}
The last deficiency, in particular, is exacerbated in battery electric shipping. The set of alternative ship types (and the associated binary variables) becomes large because every ship type and size option is duplicated for every battery capacity option.  

In this paper, we introduce a new modeling approach for ESFSMP. We model the ship sizing and speed-dependent energy consumption decisions with continuous variables and define constraints that capture their nonlinear relationships. Our main contribution is to show that these nonlinearities have a specific, favorable convex form. The resulting convex mixed-integer nonlinear programming problems can be solved in a strong sense, with a guaranteed global optimality. Our problem formulation achieves a small combinatorial search space because we only need integer-valued variables for service frequency and fleet size. Computational experiments
show that our formulation achieves 15.1\% lower fleet cost on average and reduces the solve time by at least tenfold
compared to the linear formulation.

We begin by describing the strategic planning problem in §\ref{sec:problem_description}. The next section §\ref{sec:network_model} presents the mathematical formulation of the problem using a simple base ship sizing model that scales a reference ship. The base ship design model is extended in §\ref{sec:vessel_model} with a detailed physics-based model capturing the key naval-architectural principles. We discuss the results of computational experiments in §\ref{sec:comp_experiments} and present a case study of battery-electric ro-pax ferry service design in the northern Baltic Sea in §\ref{sec:example_problem}. The reader can access the code and all the data needed to reproduce the results reported in this paper at \url{https://users.aalto.fi/~aritari/convex_ESFSMP.html}.
\section{Related work}
\subsection{Optimal fleet planning in liner shipping}
The \textit{maritime fleet size and mix problem} (MFSMP) deals with deciding the size, type, and number of ships to use to meet a given cargo demand. The objective is to minimize the total cost of purchasing and operating the fleet. MFSMPs have been extensively studied in the conventional fossil fuel-powered fleet setting \citep{CMF+07}. We limit our focus to the liner shipping segment, although MFSMP is relevant in industrial and tramp shipping as well. The variants of MFSMP are distinguished by the way ships are acquired, the evolution of the fleet over time, and the types of operating decisions that the fleet design is based on \citep{PFH14}. 

Ships can be built to exact specifications, bought in the second-hand market, or chartered. Liner container shipping companies typically operate mostly chartered ships \citep{MW11}. The MFSMP in this case consists of a short-term planning task that selects ships to be chartered from a pool of available ships. \cite{Fag99} models chartering using binary decision variables and \cite{MW11} formulates distinct fleet composition scenarios. In contrast to liner companies, ferry service providers operate a fleet of mostly owned ships. The MFSMP decision problem involves designing specifications for newbuildings to be purchased. Nevertheless, the literature is limited to purchases made from a pool of ship sizes fixed in advance \citep{LL04,SUN+05}. In contrast to the preceding studies, in this paper, we model the newbuilding characteristics using continuous variables.

In MFSMP, the fleet size and mix decisions need to be based on optimal deployment of the available ships to routes. Here, deployment refers to the number of round-trips a ship completes on a route during the planning horizon. Increasing ship speed reduces round-trip completion time but also increases energy consumption and variable operation cost. \cite{Fag99} associates a fixed design speed with each ship type. Other authors of the previously discussed works make the same assumption. Our approach in this paper is to model the nonlinear relationships between speed, transit time, and energy consumption explicitly. This approach is aligned with \cite{ZWW21} and \cite{RSH+21} from the ship routing and speed optimization literature.

\cite{HLN+22} were first to formulate MFSMP using battery electric ships (called \textit{electric ship fleet size and mix problem} (ESFSMP) in the following). Their model introduced constraints for battery capacity and charging time. The relationship between speed and energy consumption is also accounted for by linear interpolation between discrete levels as in \cite{AFH15} and \cite{FGR+15}. While \cite{HLN+22} assumed a single fixed route, the follow-up work \cite{HFN+24} introduced routing as an additional operational decision. \cite{CFJ+25} formulate an extension to multiperiod battery electric fleet planning. All the preceding battery electric ESFSMP studies use binary decision variables to select ships from a pool of available fixed ship types. However, the pool now needs to include all possible battery capacity options as well. The resulting large number of binary variables increases computational burden. In contrast, battery capacity is just another continuous variable in the nonlinear ship sizing model we formulate in this paper. 
\subsection{Battery-electric ship sizing}
Recently, many authors have investigated the limits of the economic range of battery electric ships \citep{KPP22, MPW+25}. These works estimate the battery-electric ship propulsion power by scaling the power of a reference combustion engine ship. The scaling law, known as the \textit{admiralty formula}, gives the electric ship propulsion power at speed $v$ and displacement $\nabla$ as $P_\text{prop} = P_\text{MCR} \left( \nabla/\nabla_\text{ref} \right)^{2/3} \left( v/v_\text{max} \right)^3$, where $P_\text{MCR}$ is the combustion engine power output at design speed $v_\text{max}$ and displacement $\nabla_\text{ref}$. If $v$ and $\nabla$ were assigned as decision variables in ESFSMP, the admiralty formula would give a simple nonlinear continuous model for ship sizing. However, ESFSMP needs to model size variation of at least an order of magnitude in typical applications. The formula's accuracy is insufficient when large deviations from the reference displacement can occur.

High-fidelity physics-based ship sizing model overcomes the limitations of the simple scaling law. This approach makes use of stand-alone numerical simulation tools for hydrodynamics, hydrostatics, structural analysis, and other disciplines. An integration platform defines a parametric hull form and passes the geometry data to all the connected tools \citep{Pap10}. Many authors have reported successful implementations, including damage stability-oriented design of ferries \citep{MPC+18}, tanker design with a focus on environmental performance \citep{PZB+10}, and a container ship design \citep{PBT+18}. Multidisciplinary ship design optimization using stand-alone simulation tools requires an iterative scheme with a large number of calls to the simulation routines. The computational cost becomes high, with problems taking hours or days to resolve.

Computationally expensive fluid dynamics, finite element method, and other high-fidelity simulation routines are counterproductive in ESFSMP ship sizing. Lower fidelity models of the key naval-architectural principles can provide reasonable approximations in exchange for a much lower computational burden. Recently, \cite{RMT23} demonstrated that nonlinear optimization problem formulations for high-level conceptual ship design have a special \textit{log-convex} structure. A log-convex problem becomes convex under nonlinear changes of variables. In this paper, we make use of this special structure and the efficient and reliable solution methods for convex optimization. 

\subsection{Modeling with log-convex functions}
The theoretical observation that log-transformation achieves lossless convexification of some nonconvex problems dates back to 60s to \cite{DPZ67}. Reliable and efficient numerical methods for solving log-convex problems appeared three decades later \citep{NN94}. Modeling languages for convex optimization introduced automatic variable transformations \citep{DB16}. \cite{ADB19} formulate a ruleset for constructing complex log-convex functions in a disciplined way from elementary atomic functions.

Many authors have discovered the favorable log-convex structure in important nonlinear engineering problems, including power control in wireless networks \citep{CTP+07}, gate sizing in large-scale digital circuits \citep{BKP+05}, and aircraft design \citep{HA14}. Routing cargo over a maritime transportation network can be considered analogous to transmitting over a wireless multihop network in \cite{CTP+07} or driving currents through a circuit of logic gates in \cite{BKP+05}. In the latter work, convex modeling of the maximum delay along any path through the circuit applies to the problem of determining the port operation (e.g., charging, onboarding, and loading) that takes the longest.
Similarly, the global stress response of a ship hull can be calculated with the convex beam bending model developed for wingbox structural strength analysis in \cite{HA14} 

The authors' previous work was evidently the first to illustrate the log-convex structure of ship design problems \citep{RMT23}. This work modeled a single battery-electric bulk carrier with a fixed cargo capacity and speed profile. In the present paper, we extend all the ship design subdomain models in \cite{RMT23} and embed them in the ESFSMP. 
\subsection{Positioning of the present work}
In light of the above literature discussion, our log-convex ESFSMP is categorized as a maritime fleet size and mix problem in liner shipping, where the fleet design is based on deployment decisions but not on ship routing. The linear ESFSMP studied in \citep{HLN+22} is most similar to our work. While \citep{HLN+22} artificially limits newbuilding purchase decisions to a small set of predefined types, our model allows unbounded specification of newbuildings to be purchased. Our model designs fleets with a lower total cost because the optimal ship designs are likely to lie outside of the predetermined types.

Although we introduce nonlinearities, we show that they can be reformulated in a convex form with no sacrifice in fidelity. Hence, the continuous relaxation of the mixed-integer problem is also convex. This formulation retains the key advantage of mixed-integer linear programming applied in previous work on ESFSMP: problems can be provably solved to global optimality (within a small tolerance) using a free off-the-shelf solver that requires no initial guesses or parameter tuning.

We can also mention a shortcoming of our approach. Our method does not support the encoding of logical constraints using binary variables. This restriction rules out joint optimization of fleet design and ship routing. Therefore, our method is aligned with MFSMP studies that use deployment operating decisions only \citep{MW11, HLN+22, Pre95}.
In most short-sea liner shipping segments, the routes are well established, or there are only a few viable alternatives. In this case, we can formulate an instance of our problem for every alternative route. In contrast, routing plays a larger role in industrial and tramp shipping, so our method is not well-suited to those segments.

\section{Problem description} \label{sec:problem_description}
\begin{figure*}
\begin{center}
\includegraphics[width=0.9\textwidth]{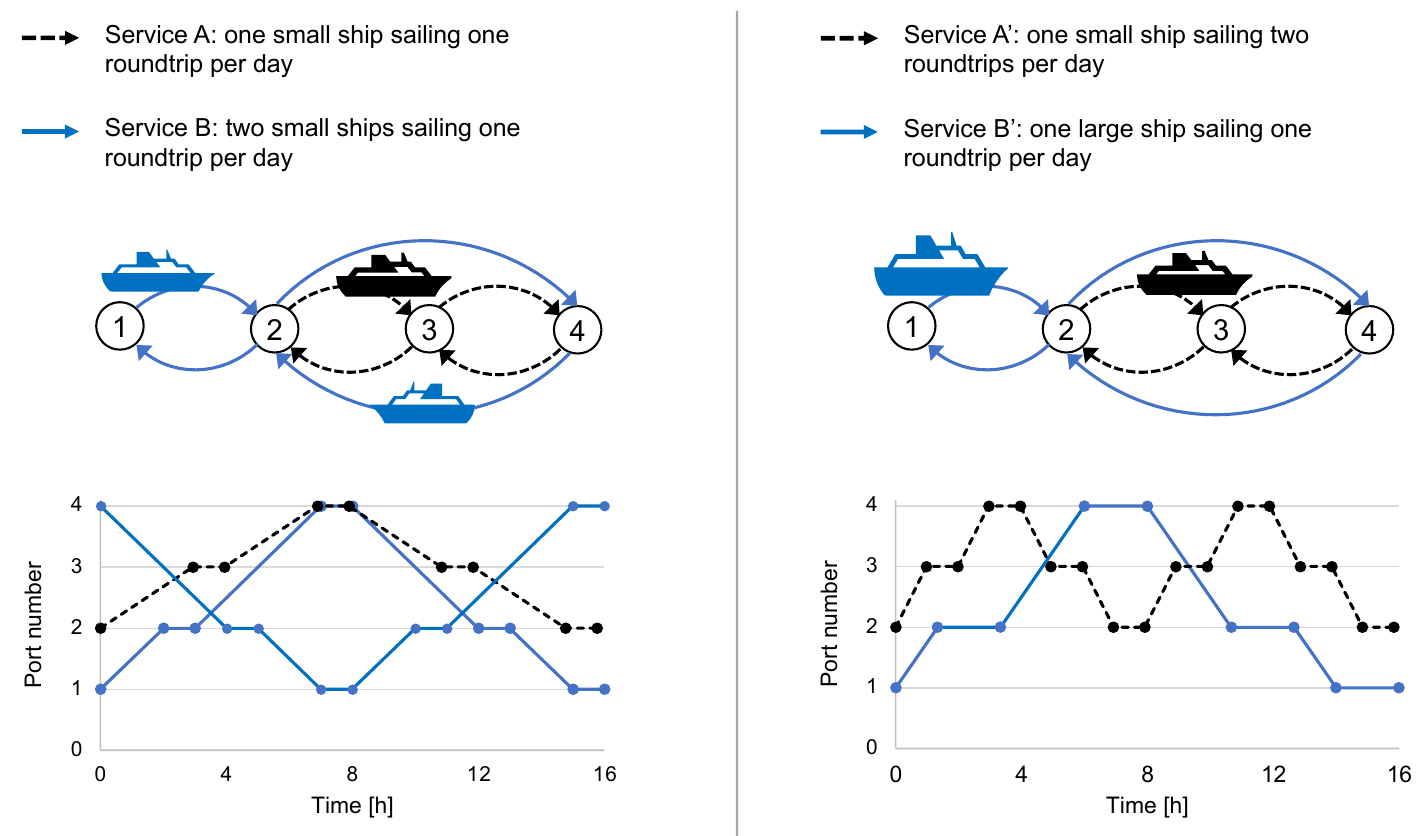} 
\caption{Two alternative fleet compositions and sailing schedules for a given route plan.} 
\label{fig:ex_routes}
\end{center}
\end{figure*}
We consider a liner shipping company operating a group of shipping routes in a given geographical area (e.g., the Baltic Sea or the Mediterranean Sea). A route consists of a sequence of ports visited in a regular cyclic schedule. One or more identical ships sailing on a route constitute a \textit{service} and all the routes and their ports form a \textit{network}. A ship can transport cargo between any origin-destination port pair in its route. We assume mixed-tonnage concept ships that can transport different types of cargoes, e.g., ro-ro cargo and passengers, as in the roll-on/roll-off passenger ferry segment.

The ports in the network are linearly ordered. Each service sails a route with a fixed origin and destination port. A route starts from a departure port and proceeds outbound to ports with higher and higher numbers until it reaches the destination port. There, the route turns around and travels inbound to ports with lower and lower numbers until it reaches the original departure port. A route does not need to call at all the ports between the end ports. However, the port calls are the same on the route's outbound and inbound directions. This is called a \textit{pendulum} service \citep{CHP+20}.

The long-term \textit{strategic} decisions are the number of ships to deploy on each service, the size of each ship, and the power rating of the onshore power supply in each port. Ship sizing includes the capacity for each type of cargo as well as battery energy capacity.  
The strategic decisions are coupled with short-term tactical and operational decisions. The frequency of each service is a tactical decision that defines the outline of the published weekly schedule. Cargo flows and sailing leg speeds are operational decisions. Energy consumption is speed-dependent, and cargo handling time is cargo flow-dependent. Port turnaround times are affected by all the operational decisions because cargo offloading, loading, and battery charging take place in parallel while the ship is docked. 

The demands of cargoes between all pairs of ports are known quantities. The demand restricts the total amount of cargo transported by all the services during the planning horizon. We also require that the onboard cargo on any sailing leg does not exceed the ship's capacity. Meanwhile, all the services must ship a minimum total weekly cargo quantity over the shipping network. This requirement is lower than the total demand, meaning that some demand is rejected. 

A round-trip is a single outbound-inbound voyage. Its duration depends on sailing leg speeds and port turnaround times. A ship cannot turn around any quicker than the port operation that takes the longest to complete. The length of the planning horizon is the time limit under which a ship must complete all its round-trips. Finally, the energy consumption on any sailing leg cannot exceed the ship's battery capacity.

The objective is to find a plan that minimizes the ship and port operators’ total cost. On the shore side, the investment cost scales with the installed power supply size. Acquiring and staffing the ships incurs a fixed cost. Energy and port charges are variable costs. They are functions of speed, ship size, and service frequency.

Figure \ref{fig:ex_routes} depicts a shipping network with four ports, two routes, and two alternative pendulum services. These examples illustrate the tradeoffs and interdependencies of strategic, tactical, and operational decisions. First, the two services are interconnected because both services can ship cargo outbound and inbound between ports 2 and 4. Second, the large ship in service B' consumes less energy per cargo ton-mile than the small ship in service B. However, the large ship has a longer turnaround time at port and a worse utilization rate than the small ship. Third, both services A and A' incur a fixed cost for one ship, but A' sails at twice the speed with twofold daily cargo capacity. This comes at a cost of at least four times greater energy consumption because resistance increases nonlinearly with speed. 
\section{Mathematical formulation} \label{sec:network_model}
\subsection{Indexing sets} \label{sec:indexing_sets}
In the mathematical description of the problem, let $\mathcal{S}$ be the set of services indexed by $s$. Ships operating a service are identical, so we also index ships by $s$. Moreover, let $\mathcal{F}$ be the set of valid frequencies indexed by $f$ and $\mathcal{C}$ be the set of cargoes indexed by $c$. 

\begin{figure}
\begin{center}
\includegraphics[width=0.5\textwidth]{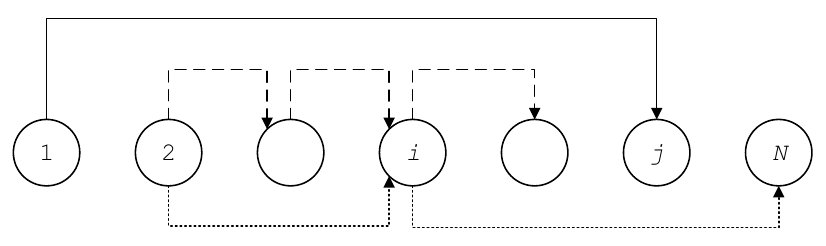} 
\caption{A route is a sequence of sailing legs between linearly ordered ports from 1 to $N$ and indexed by $i$ and $j$. The figure illustrates a shipping network consisting of three routes constructed from subsets of the same set of $N$ ports.} 
\label{fig:linear_network}
\end{center}
\end{figure}

The following notation representing ports and origin-destination demand arcs builds on the liner network design problem presented by \citet{CMF+07}. The ports in the network are linearly ordered and numbered from $1$ to $N$. Each service has two end ports $o_s$ and $d_s$, where $1 \leq o_s < d_s \leq N$ (Figure \ref{fig:linear_network}).

Let $\mathcal{N}$ be a set of the ports $1$ to $N$ in a linear order, indexed by $i,j,i'$ and $j'$.
The origin port of service $s \in \mathcal{S}$ is $o_s$, and the destination port is $d_s$.
The set $\mathcal{N}_s \subseteq \mathcal{N}$ denotes the ports of service $s \in \mathcal{S}$.

The set of ports in service $s$ that precede $i \in \mathcal{N}_s$, and include $i$, is
\begin{equation*}
    \mathcal{N}^{-}_{is}=\{i': i'\in \mathcal{N}_s, o_s \leq i' \leq i\}\subseteq \mathcal{N}_s
\end{equation*} and the ports that follow $i$ is the set
\begin{equation*}
    \mathcal{N}^{+}_{is}=\{i': i'\in \mathcal{N}_s, i \leq i' \leq d_s\}\subseteq \mathcal{N}_s.
\end{equation*}
When a ship on service $s$ departs from a port $i\in \mathcal{N}_s\backslash d_s$ in the outbound leg, its destination port is the smallest element of the set $\mathcal{N}_{is}^+\backslash \{i\}$.
Let $\mathcal{L}_s$ denote the set of outbound and inbound sailing legs in service $s$ according to 
\begin{equation*}
    \mathcal{L}_s=\left\{ (i,j),(j,i):i\in \mathcal{N}_s \backslash d_s, j=\min \mathcal{N}_{is}^+ \backslash \left\{i\right\}  \right\}.
\end{equation*}

The set of origin-destination demand arcs is
\begin{equation*}
    \mathcal{A} = \{ (i, j) : i \in \mathcal{N}, j\in \mathcal{N}, i \neq j \}.
\end{equation*}
The subset of origin-destination demands supplied by service $s$ is
\begin{equation*}
    \mathcal{A}_s = \{ (i, j) : i \in \mathcal{N}_s, j\in \mathcal{N}_s, i \neq j \} \subseteq \mathcal{A}
\end{equation*}
and the set of services supplying demand for the arc $(i, j)$ is $s:(i,j) \in \mathcal{A}_s$. A ship on service $s$ departing from port $i \in \mathcal{N}_s$ can supply any demand originating from ports $ i' \in \mathcal{N}^{-}_{is}$ and destined to ports $j \in \mathcal{N}^{+}_{is} \backslash \{i\}$.
\subsection{Parameters}
The cost of installing one unit of charger power is $C^\text{cha}$, and the cost of charging one unit of electrical energy is  $C^\text{el}$. The onshore power installed in port $i$ is upper bounded to $P^\text{cha,max}_i$. Port charge is levied in proportion to the ship's gross tonnage and denoted as $C^\text{port}$. The demand and value of a unit of cargo $c$ shipped from port $i$ to $j$ are $q^\text{dem}_{ijc}$ and $\alpha_{ijc}$. The constant $U^\text{min}$ is the total required network transportation utility. The distance of sailing directly from port $i$ to $j$ is $l^\text{leg}_{ij}$. It takes $t^\text{unit}_c$ units of time to load a unit of cargo $c$. A ship in service $s$ is available at most $t^\text{route}_s$ time during the planning horizon. Ships sailing on overnight cruises are available for all hours of the day, but in short shuttle services, they are docked overnight.  
\subsection{Decision variables}
We use the following decision variables: the positive integer variable $N^\text{ship}_s \in \mathbb{Z}_{\geq 1}, s\in\mathcal{S}$ gives the number of ships deployed in service $s$. Another integer variable $N^\text{trip}_s \in \mathcal{F}, s\in\mathcal{S}$ gives the total number of roundtrips completed in service $s$ by all the ships (each ship completes $N^\text{trip}_s/N^\text{ship}_s$ roundtrips). Variables associated with the sailing leg $(i,j)$ for ship $s$ are the sailing speed $v_{ijs} \in \mathbb{R}_{++}, (i,j)\in\mathcal{L}_s, s\in\mathcal{S}$ and the turnaround time at the destination port $t^\text{port}_{ijs}\in \mathbb{R}_{++},(i,j)\in\mathcal{L}_s, s\in\mathcal{S}$. A ship in service $s$ transports the quantity $ q_{ijcs} \in \mathbb{R}_{++}, (i,j)\in\mathcal{A}_s, c\in\mathcal{C}, s\in\mathcal{S}$ of cargo $c$ from port $i$ to $j$ during the planning horizon. The only variable associated with shore infrastructure is $P_i^\text{cha} \in \mathbb{R}_{++},i\in\mathcal{N}$, which is the installed charger power at port $i$. The ship design variables depend on the particular ship design model used. Here we represent them with the $p$-vector $x^\text{ship}_s \in \mathbb{R}_{++}^p, s\in\mathcal{S}$.
\subsection{Functions}\label{sec:funcs}
Some quantities are represented as functions of the ship design variables $x^\text{ship}_s$ and sailing leg speeds $v_{ijs}$. Let $E^\text{leg}_{ij}:\mathbb{R}_{++}^{k\times 1}\rightarrow\mathbb{R}_{++}$ denote ship electrical energy consumption sailing from port $i$ to $j$. This leg-dependent formulation accounts for sailing distance. We also define $C^\text{ship}:\mathbb{R}_{++}^{k}\rightarrow\mathbb{R}_{++}$ as ship construction and lifetime operation cost excluding energy. Ship capacity of cargo $c$ is $q^\text{cap}_c:\mathbb{R}_{++}^{k}\rightarrow\mathbb{R}_{++}$ and gross tonnage is $V^\text{GT}:\mathbb{R}_{++}^{k}\rightarrow\mathbb{R}_{++}$.
\subsection{Model formulation} \label{sec:logcvx_formulation}
\begin{table*}
\small
\begin{center}
\caption{Model notation}\label{tb:notation}
\begin{tabular}{llll} \hline
Symbol & Description & Symbol & Description \\ \hline
\multicolumn{2}{l}{\textit{Sets}} & \multicolumn{2}{l}{\textit{Parameters}} \\
$\mathcal{A}$ & Origin-destination demand arcs & $\alpha_{ijc}$ & Utility of a unit of cargo $c$ shipped \\
$\mathcal{A}_s$ & Origin-destination demand arcs of service $s$ & & from port $i$ to $j$ \\
$\mathcal{C}$ & Cargoes & $C^\text{cha}$ & Charger installation cost \\
$\mathcal{F}$ & Service frequencies & $C^\text{port}$ & Port charge \\
$\mathcal{L}_s$ & Outbound sailing legs in service $s$ & $C^\text{el}$ & Electricity price \\
$\mathcal{N}$ & Ports & $d_s$ & Outbound destination port of service $s$ \\
$\mathcal{N}_s$ & Ports visited by service $s$ & $q^\text{dem}_{ijc}$ & Demand of cargo $c$ from port $i$ to $j$ \\
$\mathcal{N}_{is}^{–}$ & Ports in service $s$ that precede port $i$ & $l^\text{leg}_{ij}$ & Sailing distance from port $i$ to $j$ \\
$\mathcal{N}_{is}^{+}$ & Ports in service $s$ that follow port $i$ & $o_s$ & Outbound departure port of service $s$  \\
$\mathcal{S}$ & Services and ships & $P^\text{cha,max}_i$ & Permissible charger power at port $i$ \\
$\mathcal{X}_s$ & Feasible values of ship $s$ design variables & $t^\text{route}_s$ &  Availability of service $s$ \\
\multicolumn{2}{l}{\textit{Decision variables}} & $t^\text{unit}_c$ & Loading duration of a unit of cargo $c$ \\
$ q_{ijcs}$ & Quantity of cargo $c$ transported from port $i$ & $U^\text{min}$ & Required transportation service utility \\
& to $j$ by ship $s$ during the planning horizon & \multicolumn{2}{l}{\textit{Functions}} \\
$N^\text{trip}_s$ & Number of roundtrips
completed in service $s$ & $C^\text{ship}$ & Ship acquisition and lifetime operation cost  \\
$N^\text{ship}_s $ & Number of ships deployed in service $s$ & $E^\text{leg}_{ij}$ & Ship electrical energy consumption  \\
$P_i^\text{cha}$ & Charger power at port $i$ &   & sailing from port $i$ to $j$  \\ 
$t^\text{port}_{ijs}$ & Turnaround time the destination port & $q^\text{cap}_c$  & Ship capacity for cargo $c$   \\ 
& of leg $(i,j)$ of ship $s$ & $V^\text{GT}$ & Ship gross tonnage \\ 
$v_{ijs}$ & Sailing speed in leg $(i, j)$ of ship $s$ & &  \\ 
$x^\text{ship}_s $ & Vector of ship $s$ design variables & &  \\
\hline
\end{tabular}
\end{center}
\end{table*}
The model formulation is as follows:

\phantom{Some text here that flows to the next line a bit more}
minimize
\begin{equation}
    \begin{aligned}
        & \sum_{s \in \mathcal{S}} \Biggl[N^\text{ship}_s C^\text{ship}(x^\text{ship}_s)
        + \sum_{(i, j) \in \mathcal{L}_s} N^\text{trip}_s\Bigl( C_j^\text{port} V^\text{GT}(x^\text{ship}_s) \\
        &+ C_j^\text{el} E^\text{leg}_{ij}(x^\text{ship}_s, v_{ijs}) \Bigr)\Biggr]
        + \sum_{i' \in \mathcal{N}}  C^\text{cha}P_{i'}^\text{cha},   \label{eq:obj}
    \end{aligned}
\end{equation}
subject to
\begin{align}
& \sum_{i'\in \mathcal{N}^{-}_{is}} \sum_{j'\in \mathcal{N}^{+}_{is} \backslash \{i\}} q_{i'j'cs} \leq q_{c}^\text{cap}(x^\text{ship}_s),
    \begin{array}{ll}
         & c\in \mathcal{C}, \, s\in \mathcal{S}, \\
         & i\in \mathcal{N}_s \backslash \{d_s\},
    \end{array} \label{eq:constr1} \\
& \sum_{i'\in \mathcal{N}^{-}_{is}} \sum_{j'\in \mathcal{N}^{+}_{is} \backslash \{i\}} q_{j'i'cs} \leq q_{c}^\text{cap}(x^\text{ship}_s), 
\begin{array}{ll}
         & c\in \mathcal{C}, \, s\in \mathcal{S}, \\
         & i\in \mathcal{N}_s \backslash \{d_s\},
    \end{array} \label{eq:constr2} \\
& \sum_{s:(i,j) \in \mathcal{A}_s} N^\text{trip}_s q_{ijcs} \leq q_{ijc}^\text{dem},\quad c\in \mathcal{C},\ (i,j) \in \mathcal{A}, \label{eq:constr3} \\
& \sum_{(i,j) \in \mathcal{L}_s} N^\text{trip}_s (t_{ijs}^\text{transit} + t_{ijs}^\text{port}) \leq N^\text{ship}_st_s^\text{route}, \quad s\in \mathcal{S}, \label{eq:constr4} \\
& \max \{t_{ijs}^\text{cargo}, t_{ijs}^\text{cha} \} \leq t_{ijs}^\text{port}, 
     \quad (i,j) \in \mathcal{L}_s,\, s\in \mathcal{S}, \label{eq:time3} \\
& \sum_{ (i,j) \in \mathcal{A} } \sum_{s:(i,j) \in \mathcal{A}_s} \sum_{c\in \mathcal{C}} \alpha_{ijc} \log (N^\text{trip}_sq_{ijcs})\geq U^\text{min}, \label{eq:constr5} \\
&\quad x^\text{ship}_s \in \mathcal{X}_s,\, s \in \mathcal{S}, \quad P^\text{cha}_i \leq P^\text{cha,max}_i, \, i\in\mathcal{N}, \label{eq:constr6}
\end{align}
where $t_{ijs}^\text{transit},\ldots,t_{ijs}^\text{cha}$ are to be substituted with the expressions
\begin{align}
    & t_{ijs}^\text{transit}=\frac{l_{ij}^\text{leg}}{v_{ijs}}, \label{eq:time2} \\
    & t_{ijs}^\text{cargo} = \max\{ t_{ijsc}^\text{load}+t_{ijsc}^\text{unload}: c \in \mathcal{C}\}, \label{eq:time7} \\
    & t_{ijsc}^\text{unload} = t_{c}^\text{unit} \sum_{i'\in \mathcal{N}^{-}_{js} \backslash \{j\}} q_{i'jcs}, \label{eq:time4} \\
    & t_{ijsc}^\text{load} = \left\{ \begin{array}{lr} 
    t_{c}^\text{unit} \sum_{j'\in \mathcal{N}_{s} \backslash \{j\}} q_{jj'cs} & \text{if } j = d_s \\ 
    t_{c}^\text{unit} \sum_{j'\in \mathcal{N}^{+}_{js} \backslash \{j\}} q_{jj'cs} & \text{otherwise,}
         \end{array} \right. \label{eq:time5} \\
    & t_{ijs}^\text{cha}=\frac{E^\text{leg}_{ij}(x^\text{ship}_s, v_{ijs})}{P^\text{cha}_j} \label{eq:time6}.
\end{align}
The objective function \eqref{eq:obj} minimizes the total annualized cost. The four terms are the ship purchase and running cost, port charges, electricity cost, and shore charger installation cost. All these costs are scaled to one year.

The constraints are as follows:
\begin{itemize}
    \item The capacity constraints \eqref{eq:constr1} ensure that the ship capacity is not exceeded for any cargo type on the outbound legs. For the inbound legs \eqref{eq:constr2}, the capacity constraints take the same form, with the exception that the departure and arrival ports appear in a reverse order.
    \item The demand constraints \eqref{eq:constr3} state that the demand of an origin-destination port pair is the upper bound on the total quantity of cargo transported by all the services that include both ports in their route.
    \item The upper bound on voyage duration is expressed in \eqref{eq:constr4}. The duration of a one-way voyage is the sum of times on sea \eqref{eq:time2} and port. Cargo unloading \eqref{eq:time4}, loading \eqref{eq:time5}, and charging \eqref{eq:time6} take place in parallel, and the operation that takes the longest determines the turnaround time \eqref{eq:time3}. 
    \item We require that the transportation network provide a minimum service level \eqref{eq:constr5}. Here, we use log-utility, which is a particular concave function with diminishing incremental value of cargo volume shipped in each origin-destination pair and cargo type.
\end{itemize}
\subsection{Base ship design model}
We now describe a particular simple ship sizing model for the ESFSMP defined in the previous section. This model uses a single sizing variable that scales the displacement of a reference ship. The accuracy is acceptable only for ship sizes within a small neighborhood of the reference ship. 

The vector of design variables is $x^\text{ship}_s=(\nabla_s,E^\text{max}_s)$, where $\nabla$ is the displacement volume and $E^\text{max}_s$ is the battery capacity. Energy consumption on leg $(i,j)$ using speed $v_{ijs}$ is
\begin{equation*}
    E_{ijs}^\text{leg} = t^\text{transit}_{ijs} P^\text{ref} \left( \frac{\nabla}{\nabla^\text{ref}} \right)^{2/3} \left( \frac{v_{ijs}}{v^\text{ref}} \right)^3,
\end{equation*}
where $P^\text{ref}$ is the propulsion power of a reference ship at $v^\text{ref}$ and displacement $\nabla^\text{ref}$. The feasible set is
\begin{equation} \label{eq:admiralty_feasibl_set}
    \mathcal{X}_s=\{x^\text{ship}_s:E_{ijs}^\text{leg}\leq E^\text{max}_s,(i,j)\in\mathcal{L}_s\},\quad s\in\mathcal{S}.
\end{equation}
The functions in §\ref{sec:funcs} take a monomial form, e.g., $q^\text{cap}_{c}(\nabla_s)=\Psi_c^{(1)}(\nabla_s)^{\Psi_c^{(2)}}$ with parameters $\Psi_c^{(1)},\Psi_c^{(2)}$.
\subsection{Convex formulation} \label{ref:cvx_formulation}
The ESFSMP described by \eqref{eq:obj} to \eqref{eq:admiralty_feasibl_set} is a finite-dimensional, nonlinear, mixed-integer programming problem. It has, however, a special structure. The objective and the constraints are defined by log-convex functions. Hence, the problem is a mixed-integer log-convex optimization problem. As a result, it is transformed to a \textit{convex} mixed-integer nonlinear programming problem using simple changes of variables. This transformation is exact, i.e., it does not involve any approximations. The transformed problem can be solved in a strong sense, with a guaranteed global optimality. For background on formulation techniques and solution methods for convex mixed-integer nonlinear programming, see \cite{KBL+19}.

In the following, we describe in detail the log-convex structure of the ESFSMP. Consider a particular log-convex function, called a \text{posynomial}, defined as
\begin{equation} \label{eq:posy_form}
    f(x)=\sum_{k=1}^Kc_kx_1^{a_{1k}}x_2^{a_{2k}}\cdots x_n^{a_{nk}},
\end{equation}
where $x\in\mathbb{R}^n_{++}$ is the $n$-vector of variables, $a_{1k},\ldots,a_{nk}\in\mathbb{R}$ and $c_k\in\mathbb{R}_{++}$ are constants and $K$ is the number of product terms. The ESFSMP can be expressed as a problem that minimizes a posynomial over a set defined by posynomial inequalities and the restriction that some number $P\leq n$ of variables are positive integers:
\begin{equation} \label{eq:milc}
    \begin{aligned}
        &\text{minimize} \quad && f_0(x) \\
        &\text{subject to} \quad && f_m(x) \leq 1, \quad m = 1,\ldots,M, \\
        &                        && x_l \in \mathbb{Z}_{\geq 1}, \quad l=1,\ldots,P,
    \end{aligned}
\end{equation}
where $f_0,\ldots,f_P$ are posynomials. This problem becomes a convex mixed-integer nonlinear programming problem under a logarithmic change of variables, see \ref{sec:var_changes}.

To see that the ESFSMP is defined by posynomials, consider the following reformulations. First, the capacity constraint \eqref{eq:constr1}, where $q^\text{cap}_{c}(\nabla_s)=\Psi_c^{(1)}(\nabla_s)^{\Psi_c^{(2)}}$, takes the posynomial inequality form \eqref{eq:milc} by a simple algebraic manipulation:
\begin{equation*}
    \sum_{i'\in \mathcal{N}^{-}_{is}} \sum_{j'\in \mathcal{N}^{+}_{is} \backslash \{i\}} \frac{q_{i'j'cs}}{\Psi_c^{(1)}\nabla_s^{\Psi_c^{(2)}}} \leq 1,
    \begin{array}{ll}
         & c\in \mathcal{C}, \, s\in \mathcal{S}, \\
         & i\in \mathcal{N}_s \backslash \{d_s\}.
    \end{array}
\end{equation*}
The utility constraint \eqref{eq:constr5} is equivalent to
\begin{equation*}
\prod_{(i,j) \in \mathcal{A}}  \prod_{s:(i,j) \in \mathcal{A}_s} \prod_{c\in \mathcal{C}}  (N^\text{trip}_sq_{ijcs})^{\alpha_{ijc}} \geq \exp{U^\text{min}}.
\end{equation*}
The port turnaround time constraint \eqref{eq:time3} is equivalent to the inequalities
\begin{equation*}
    \begin{aligned}
        & t_{ijs}^\text{port} \geq t_{ijs}^\text{load}, \quad (i,j) \in \mathcal{L}_s,\, c\in \mathcal{C},\,  s\in \mathcal{S}, \\ 
        & t_{ijs}^\text{port} \geq t_{ijs}^\text{unload}, \quad (i,j) \in \mathcal{L}_s,\, c\in \mathcal{C},\,  s\in \mathcal{S}, \\ 
        & t_{ijs}^\text{port} \geq t_{ijs}^\text{cha}, 
     \quad (i,j) \in \mathcal{L}_s,\, s\in \mathcal{S}.
    \end{aligned}
\end{equation*}
\section{Physics-based ship design extension} \label{sec:vessel_model}
\begin{table}
\small
\begin{center}
\caption{New decision variables for the extended ship design model. We use $s$ to index both services and ships because all the ships deployed on a service are identical.}\label{tb:vars}
\begin{tabular}{ll} \hline
Symbol & Description  \\ \hline
$L_s$ & Hull length \\
$L^\text{sup}_s$ & Superstructure length \\
$B_s$ & Hull waterline breadth \\ 
$T_s$ & Hull draft \\
$D_s$ & Hull depth  \\
$p_s^\text{td}$ & Reference plate thickness  \\ 
$h^\text{batt}_s$ & Battery room height\\
$w^\text{batt}_{ks}$ & Width of battery room $k$\\
$l^\text{batt}_{ks}$ & Length of battery room $k$ \\
$\tilde{h}^\text{batt}_s$ & Battery room floor elevation from keel \\ 
$\tilde{l}^\text{batt}_{ks}$ & Battery room $k$ front wall distance \\ 
& from forward perpendicular \\
$\tilde{h}^\text{roro}_s$ & First ro-ro deck elevation from keel \\
$E_s^\text{batt}$ & Battery capacity \\ 
$t_s^\text{cell}$ & Cell cycle life \\ 
$N_s^\text{batt}$ & Number of battery replacements \\  
\hline
$r_{ijs}^\text{Cf}$ & Auxiliary variable for modeling frictional \\
& resistance in sailing leg $(i,j)$\\
$r_{ijs}^\text{Cr,std}$ & Auxiliary variable for modeling residual \\
& resistance in sailing leg $(i,j)$ \\
$r^g_{k's}$ & Auxiliary variable for point $k'$ in the numerical \\
& integration of hull wetted area \\
$r^{g2}_{k's}$ & Auxiliary variable for point $k'$ in the numerical \\
& integration of hull steel area \\
\hline
\end{tabular}
\end{center}
\end{table}
This section describes a physics-based ship model that captures the key naval-architectural principles governing ship design. Although we focus on modeling a ro-pax ferry, the model applies to other ship types with minor modifications. 

The model is expressed as a group of posynomial inequalities, meaning we can manipulate them to the form $f(x)\leq1$, where $f$ is a posynomial function \eqref{eq:posy_form}. Hence, the model we describe is compatible with the ESFSMP in §\ref{sec:logcvx_formulation}. Table \ref{tb:vars} lists the new decision variables that constitute the elements of the vector $x^\text{ship}$. The domain of all the variables is $\mathbb{R}_{++}$. All the constraints and the associated design variables hold for each $s\in \mathcal{S}$, but we drop the index for clarity. Note that some symbols used in §\ref{ref:cvx_formulation} are redefined in this section.
\subsection{Weight and buoyancy equilibrium} \label{sec:weight_buo_balance}
The most elementary design relation states that buoyancy and weight in calm water are in balance because the ship is a free-floating body. This is expressed as
\begin{equation} \label{eq:vertical_force_balance}
    \nabla \rho^\text{sw} \geq W_L + W_D,
\end{equation}
where $\nabla$ is the displaced volume at design draft, $\rho^\text{sw}$ is the sea water density, and $W_L$ and $W_D$ are lightweight and deadweight in tons. We give a detailed breakdown of the weight in §\ref{sec:weight_breakdown} and derive a monomial expression for $\nabla$ in the following.

Surfaces defined by closed-form functions describe the submerged hull form. The origin of the coordinate system is located at the foremost point at the keel centerline in a Cartesian coordinate system $XYZ$.
The $X$-axis points to the port side, the $Y$-axis points to the direction of the stern, and the $Z$-axis points upwards.
Let $H^\text{fore}:\mathbb{R}_{+} \times \mathbb{R}_{+} \rightarrow \mathbb{R}_{+}$ and $H^\text{aft}:\mathbb{R}_{+} \times \mathbb{R}_{+} \rightarrow \mathbb{R}_{+}$ denote two functions that return the lateral offset from the centerline. The particular functions we use are
\begin{align}
        H^\text{fore}(y,z)=&\frac{B}{2} \sqrt{\frac{2y}{L}} \left( \frac{z}{T} \right)^{1/\beta}, \label{eq:hull_form1} \\
        H^\text{aft}(y,z)=&\frac{B}{2}\left( \frac{z}{T} \right)^{1/\beta(1-2y/L)}, \quad y\neq L/2, \label{eq:hull_form2}
\end{align}
where $L$ is the length overall, $T$ is the design draft and $B$ is the waterline breadth. All these dimensions are decision variables.
In \eqref{eq:hull_form1} and \eqref{eq:hull_form2}, $\beta$ is a hull form parameter with a higher value corresponding to a higher block coefficient. Using the above functions, we can describe the half-breadth hull form as a surface in $\mathbb{R}^3$ as
\begin{equation*}
    \begin{aligned}
          &\{ (H^\text{fore}(y,z),y,z):0\leq y \leq L/2, 0\leq z \leq T \} \\ \cup \, & \{ (H^\text{aft}(y-L/2,z),y,z):L/2\leq y < L, 0\leq z \leq T \}.
    \end{aligned}
\end{equation*}

\begin{figure}
\begin{center}
\includegraphics[width=0.4\textwidth]{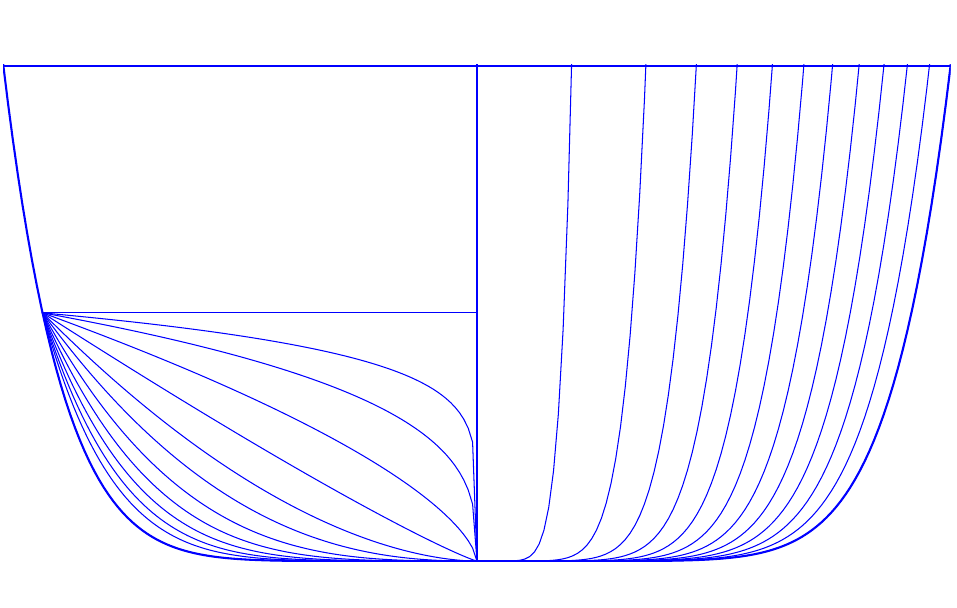} 
\caption{Line drawing of a hull defined by offset functions. The lines represent the intersection of the hull with transverse planes. Lines on the right-hand side of the vertical centerline represent the forward section, and the left-hand side the aft section.} 
\label{fig:body_plan}
\end{center}
\end{figure}

Figure \ref{fig:body_plan} illustrates the submerged hull geometry with $\beta=6$.
In this illustration, the hull is extended from the waterline to the topmost continuous deck using the geometry of the forward section.

Displacement volume, wetted surface area, hull steel weight, and buoyant force action point are functions of the variables $L, B, T$ and the form parameter $\beta$.
The hull displacement volume $\nabla$ at draft $T$ is
    \begin{equation*}
        \begin{aligned}
             \nabla & = \int_{0}^{L/2} \int_{0}^{T} (H^\text{fore}(y,z) + H^\text{aft}(y,z))\ dz\ dy \\ 
             & = \underbrace{ \left( \frac{1}{2} + \frac{\beta}{3(1+\beta)} - \frac{\log{(\beta+1)}}{2\beta} \right)}_{C_B} TLB,
        \end{aligned}
    \end{equation*}
where $C_B$ is the block coefficient.
\subsection{Cargo capacity}
We model mixed-tonnage concept ships with dedicated spaces for ro-ro cargo and passengers. In our notation, the set $\mathcal{C}$ is $\{\text{pax, roro} \}$. Here, we don't distinguish between passenger vehicles, trucks, trailers, buses, and other types of ro-ro cargo.

\begin{figure}
\begin{center}
\includegraphics[width=0.45\textwidth]{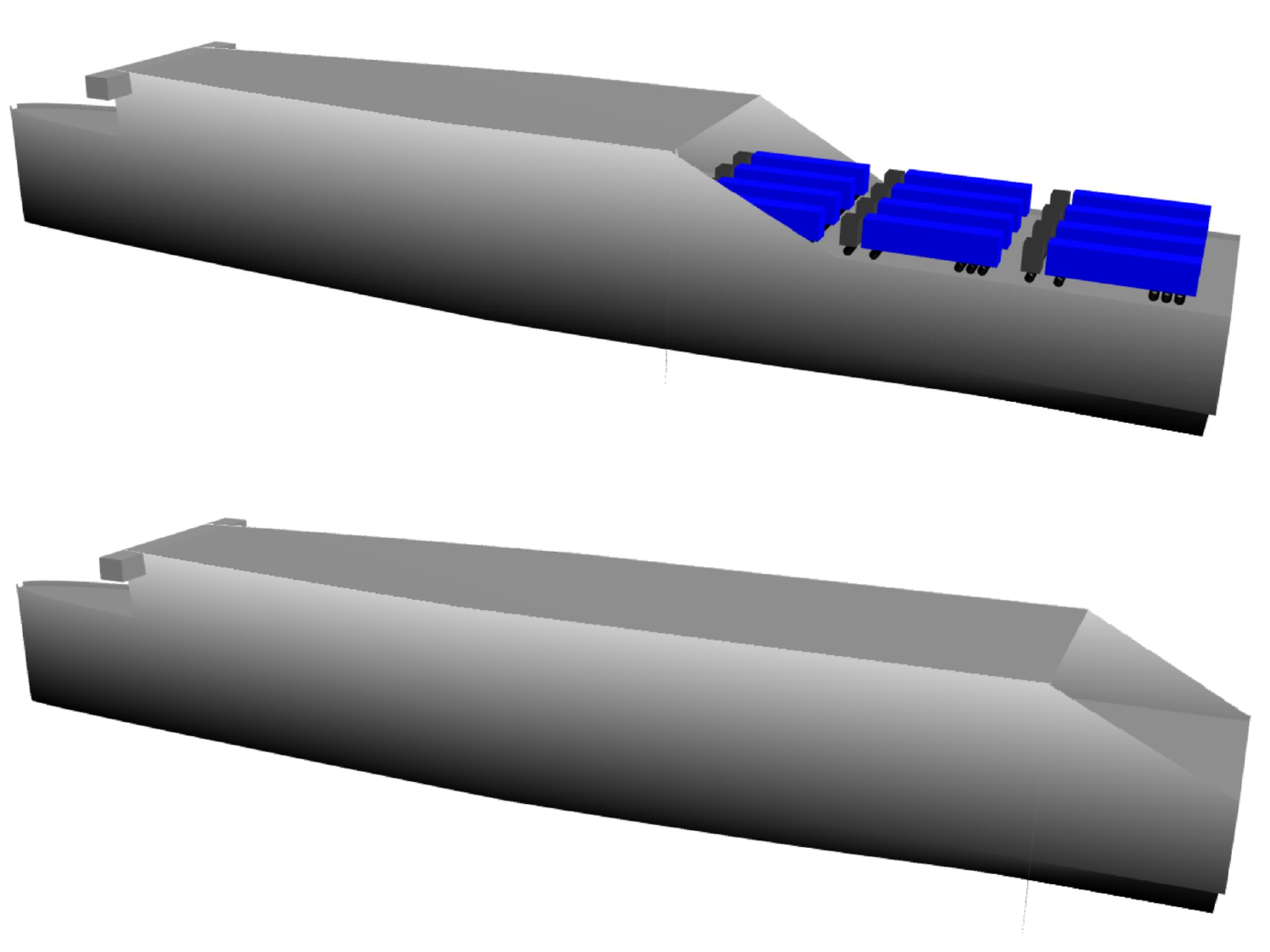}
\caption{Superstructure length as a design variable. Longer superstructure increases passenger capacity, but also electricity consumption, steel weight, construction cost, and port charges.} 
\label{fig:ropax_hull}
\end{center}
\end{figure}

Passengers' accommodation floor area and ro-ro lane capacity are functions of the hull principal dimensions, superstructure length (figure \ref{fig:ropax_hull}), and the number of superstructure decks $N^\text{sup}$. The corresponding capacities $f_{c=\text{pax}}^\text{cap}$ and $f_{c=\text{roro}}^\text{cap}$ are expressed as 
\begin{equation*}
        f_{c=\text{pax}}^\text{cap}=\frac{N^\text{sup}BL^\text{sup}}{\varphi^\text{pax}}, \quad f_{c=\text{roro}}^\text{cap}=\frac{2BL}{\varphi^\text{roro}}.
\end{equation*}
The parameters $\varphi^\text{pax}$ and $\varphi^\text{roro}$ are unit areas. 
\subsection{Upright transverse stability}
\begin{figure}
\begin{center}
\includegraphics[width=0.4\textwidth]{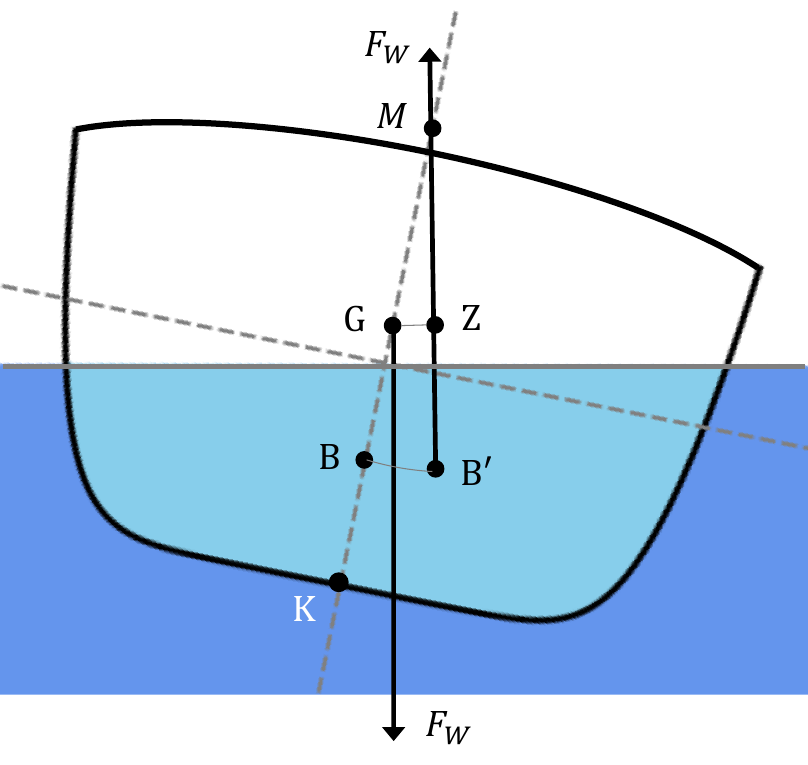}
\caption{Shifting of the center of buoyancy due to heel. The lever arm GZ between buoyancy and weight vectors must be positive by design. Under this condition, the ship reverts to upright orientation.} 
\label{fig:stability}
\end{center}
\end{figure}

Stability measure calculation relies on the notion of \textit{metacenter}. It is the position on the hull vertical centerline where the buoyancy force vector and the centerline intersect for a small angle of heel. Metacentric height $\tilde{h}_\text{GM}$ is the vertical distance from the center of gravity to the metacenter:
\begin{equation*}
    \tilde{h}_\text{GM} = \tilde{h}_\text{KB} + \tilde{h}_\text{BM} - \tilde{h}_\text{KG},
\end{equation*}
where $\tilde{h}_\text{KB}$ is the distance from keel to the center of buoyancy, $\tilde{h}_\text{BM}$ is the distance from the center of buoyancy to the metacenter and $\tilde{h}_\text{KG}$ is the distance from keel to the center of gravity (figure \ref{fig:stability}). The ship is stable when $\tilde{h}_\text{GM}>0$ by design, meaning the ship reverts to an upright orientation from a small angle of heel.

The vertical position of the centroid of the displaced water is
    \begin{equation} \label{eq:centroid_vert_pos}
        \begin{aligned}
                \tilde{h}_\text{KB} & = \frac{1}{\nabla} \int_{0}^{L/2} \int_{0}^{T} z(H^\text{fore}(y,z) + H^\text{aft}(y,z))\ dz\ dy  \\ 
                & = \left[ \frac{ \beta \hat{\beta} \tilde{\beta}/2 - \hat{\beta} \tilde{\beta} \log{\tilde{\beta}}/4 + 2\beta^2 \hat{\beta}/3 }{\beta \hat{\beta} \tilde{\beta} - \hat{\beta} \tilde{\beta} \log{\hat{\beta}} + 2\beta^2 \tilde{\beta}/3  } \right]T=\beta_\text{KB}T,
            \end{aligned}
    \end{equation}
where $\hat{\beta}=\beta+1$ and $\tilde{\beta}=2\beta+1$

The vertical distance from the center of buoyancy to the metacenter is $\tilde{h}_\text{BM} = I^\text{wp}/\nabla$, where $I^\text{wp}$ is the area moment of inertia of the waterplane area about its centroidal axis. We compute the (definite) moment integral
\begin{equation*}
    I^\text{wp} = \frac{2}{3}\int_0^{L/2} (H^\text{fore}(y,T)^3 + H^\text{aft}(y,T)^3) \, dy = \frac{7LB^3}{120},
\end{equation*}
which leads to
\begin{equation*}
    \tilde{h}_\text{BM} = \frac{7LB^3}{120C_BLBT} = \frac{7B^2}{120C_BT}.
\end{equation*}
The positive metacentric height requirement for stability follows from the quantities defined above as 
\begin{equation*}
   \beta_\text{KB}T + \frac{7B^2}{120C_BT} -\tilde{h}_\text{KG} \geq \epsilon_\text{GM},
\end{equation*}
where $\epsilon_\text{GM}$ is small and positive. Reorganizing the terms gives
\begin{equation*}
    1 + \frac{7}{120C_B\beta_\text{KB}} \left(\frac{B}{T} \right)^2 \geq \frac{\epsilon_\text{GM} + \tilde{h}_\text{KG}}{\beta_\text{KB}T}. 
\end{equation*}
This constraint is not a valid log-convex constraint because it imposes a lower bound on a posynomial. We fit a monomial to the posynomial in the relevant range of feasible $B/T$ values (roughly $[2.5, 6]$) and rewrite the constraint as
\begin{equation} \label{eq:transverse_stability}
    \delta_1 \left(\frac{B}{T} \right)^{\delta_2} \geq \frac{\epsilon_\text{GM} + \tilde{h}_\text{KG}}{\beta_\text{KB}T}. 
\end{equation}
Here $\delta_1$ and $\delta_2$ are the coefficients of the fitted monomial. The root mean square error is less than 1\%.
\subsection{Structural strength}
The structural strength model ensures that the hull withstands external bending and shear loads without failure. 
Global structural behavior is the response of the hull to loads acting at the entire length of the hull. 
At this level, the ship is modeled as a beam, called \textit{hull girder}, that is loaded with deadweight and lightweight and supported by buoyancy. 
The properties of the midship cross-section determine the strength of the hull girder. 

\begin{figure*}
\begin{center}
\includegraphics[width=0.7\textwidth]{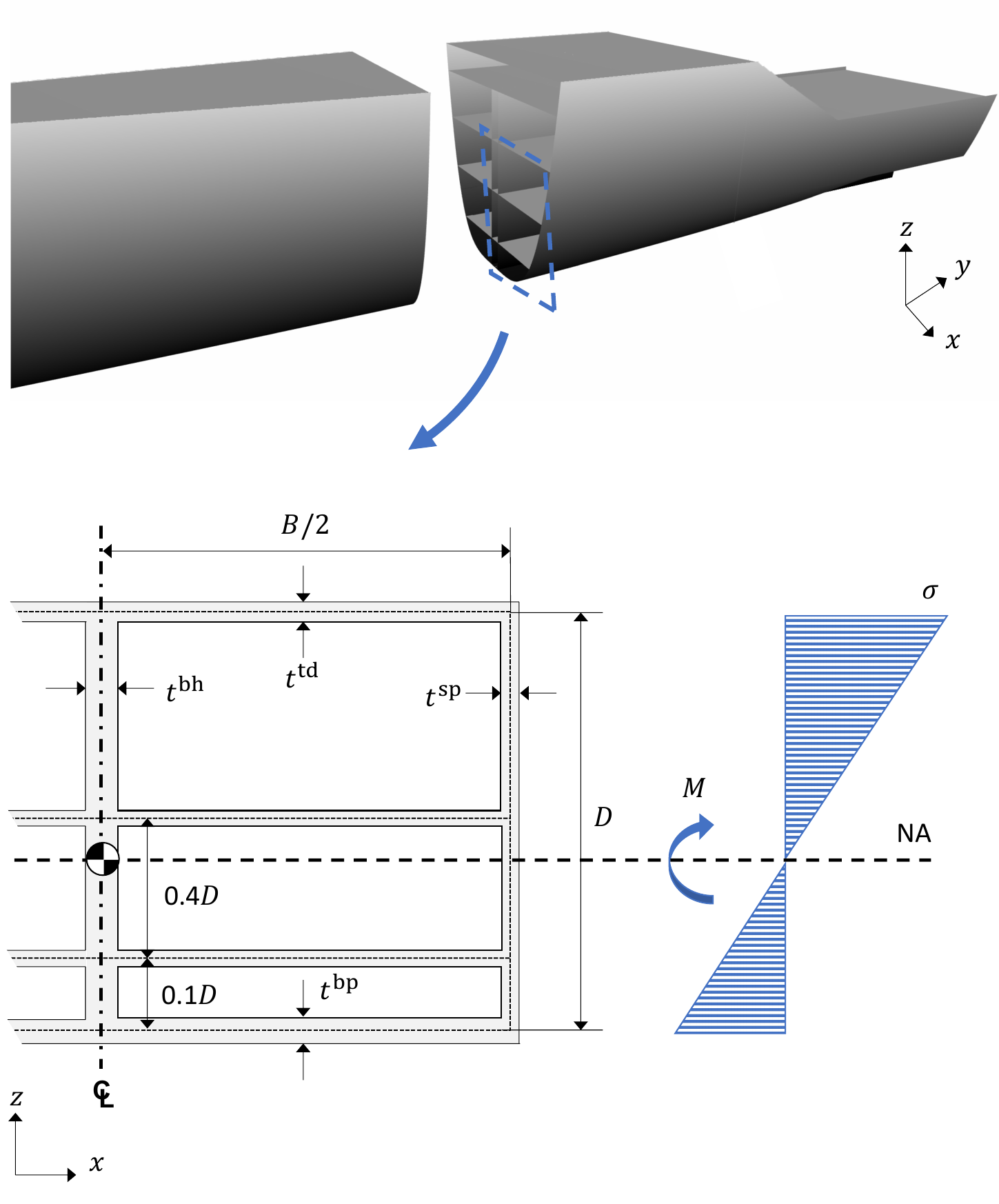}
\caption{Idealized hull girder cross-section. The normal stress due to vertical bending moment is linearly distributed across the section.} 
\label{fig:cross_section}
\end{center}
\end{figure*}

Figure \ref{fig:cross_section} depicts the hull girder cross-section topology of a ro-pax ferry. 
The girder extends from the keel to the topmost continuous deck, excluding the superstructure, which does not span the entire length of the hull. 
Only continuous longitudinal elements of the hull contribute to the structural strength. 
These include the bottom, inner bottom, side shell, longitudinal bulkheads, and the ro-ro decks. 
The elements and their relative positions are fixed, but principal dimensions and plate thicknesses are decision variables whose values determine the strength properties. 

An imbalance of buoyancy and weight along the length of the hull creates a bending moment that rises to the maximum value at amidships. We integrate twice the net force per unit length (estimated in \ref{sec:sw_loads}), which results in the approximate still water bending moment at amidships:
\begin{equation*}
    M^\text{sw} = 0.0472L^2B(C_B+0.7).
\end{equation*}
A ship on waves experiences slowly varying buoyancy forces. Classification rules give the formula for the maximum wave bending moment according to the geometry of the hull:
\begin{equation*}
    M^\text{wv}=\varphi^\text{sh}0.975L^2B(C_B+0.7),
\end{equation*}
where $\varphi^\text{sh}=-1.1$ for sagging and $ \varphi^\text{sh}=1.9C_B/(C_B+0.7)$ for hogging.
By taking the net design bending moment as $M^\text{des}=\lvert M^\text{wv}\rvert + \lvert M^\text{sw} \rvert$, we can impose the high level bending stress constraint
\begin{equation} \label{eq:bending_strength}
    \frac{M^{\text{des}}}{Z^\text{deck}} \leq \sigma^\text{perm},
\end{equation}
where $\sigma^\text{perm}$ is the permissible normal stress. The quantity $Z^\text{deck}$ is a monomial expression for the deck section modulus, which relates the applied moment to the maximum bending stress. It depends on the midship cross-section and is derived in \ref{sec:sec_modulus}.

The still water shear force at the longitudinal coordinate $y$ is obtained from the relation $V^\text{sw}(y)=dM^\text{sw}(y)/dy$. The wave-induced shear force is
\begin{equation*}
    V^\text{wv}(y)=\varphi_V(y)0.3LB(C_B+0.7),
\end{equation*}
where the value of the coefficient $\varphi_V\in[-0.92,1]$ varies along the length of the hull. Classification rules are given for both positive and negative shear forces. By taking the still water shear force, and adding to it the contribution from waves, we find the net design shear force $V^\text{des} = \lvert V^\text{sw} \lvert + \lvert V^\text{wv} \rvert$ located at $y=0.3L$. The corresponding stress limit constraint is
\begin{equation} \label{eq:shear_strength}
    \lvert q_\text{NA}^\text{sp} \rvert \frac{V^\text{des}}{p^\text{sp}} \leq \tau^\text{perm}.
\end{equation}
Here $\tau^\text{perm}$ is the permissible shear stress. The vertical shear force is mapped to the maximum shear stress response via the quantity $q_\text{NA}^\text{sp}$, which is a monomial expression for the shear flow at the neutral axis of the side plating. It is derived in \ref{sec:shear_flow}.
\subsection{Resistance}
The total ship resistance $R$ is modeled as a sum of contributions from two sources: frictional resistance $R_F$ and residual resistance $R_R$. We derive both components in the following.
The total resistance enters into the effective towing power formula in §\ref{sec:batt_model}, which determines the power drawn from the battery along with propeller and power conversion efficiency.

Frictional resistance is
    \begin{equation*}
            R_F =\frac{1}{2} C_F\rho^\text{sw}v^2A_S,
    \end{equation*}
where the coefficient $C_F$ can be estimated by the ITTC-57 method via the Reynold's number $\mathrm{Re}$ and the kinematic viscosity of water $\nu$:
    \begin{equation*}
            C_F=\frac{75}{(\log_{10}(\mathrm{Re})-2)^2}, \quad \mathrm{Re}=\frac{vL}{\nu}.
    \end{equation*}
We obtain a log-convex expression for $C_F$ by introducing the auxiliary variable $r_{C_F}$:
    \begin{equation} \label{eq:fric_res}
            C_F=\frac{75}{r_{C_F}^2}, \quad r_{C_F} + 2 \leq \log_{10}(\mathrm{Re}).
    \end{equation}
A numerical integration scheme for calculating the wetted area $A_S$ is described in \ref{sec:wetted_area}.
The residual resistance is given by (note $BT$ as the reference area instead of the wetted surface):
    \begin{equation*}
            R_R = C_R\frac{\rho}{2}v^2BT.
    \end{equation*}
The empirical method by \cite{Hol98} gives the coefficient $C_R$ for a ship floating on even keel as
    \begin{equation} \label{eq:hllb}
            C_R = C_{R}^\text{std}C_{R}^\text{Fncrit}k_L\left( \frac{T}{B} \right)^{\kappa_1} \left( \frac{B}{L} \right)^{\kappa_2} \left( \frac{D_P}{T} \right)^{\kappa_3}\left(N^\mathrm{rud}\right)^{\kappa_4},
    \end{equation}
where $N^\mathrm{rud}$ is the number of rudders and $\kappa_1,\ldots,\kappa_4$ are constants that depend on the hull type (single or twin-screw configuration). The length dependency factor $k_L$ is the monomial $k_L = \psi_1 L^{\psi_2}$ with constants $\psi_1$ and $\psi_2$.

The standard coefficient $C_{R}^\text{std}$ is expressed as a function of the Froude number $\mathrm{Fr}=v/\sqrt{gL}$ with constants $\omega_1,\ldots,\omega_9$ as
    \begin{equation} \label{eq:crstand}
            \begin{aligned}
                C_{R}^\text{std} &= \omega_1 + \omega_2\mathrm{Fr} + \omega_3\mathrm{Fr}^2 + C_B(\omega_4 + \omega_5\mathrm{Fr} + \omega_6\mathrm{Fr}^2) \\
                & + C_B^2 (\omega_7 + \omega_8\mathrm{Fr} + \omega_9\mathrm{Fr}^2).
            \end{aligned}
    \end{equation}

Since the sum $\omega_2 + C_B\omega_5 + C_B^2\omega_8$ is negative according to \cite{Hol98}, \eqref{eq:crstand} cannot be expressed as a log-convex constraint directly.
However, we can fit a softmax-affine function \citep{HKA16} to log-transformed data generated with \eqref{eq:crstand} in the relevant $\mathrm{Fr}$ range. 
A fitting function with two exponential terms achieves a sufficiently good fit. The root mean square error is less than 0.1\%.

The coefficient $C_R^\text{Frcrit}$ is expressed as
    \begin{equation} \label{eq:crcrit_orig}
            C_{R}^\mathrm{Frcrit} = \max \left\{ 1.0, \left( \frac{\mathrm{Fr}}{\mathrm{Fr}^{\mathrm{crit}}} \right)^{\left( \frac{\mathrm{Fr}}{\mathrm{Fr}^{\mathrm{crit}}} \right)}  \right\},
    \end{equation}
where $\mathrm{Fr}^{\mathrm{crit}} = \theta_1 + \theta_2 C_B + \theta_3 C_B^2$ with constants $\theta_1,\ldots,\theta_3$.
Since $C_B$ is constant for a given hull form parameter $\beta$, $\mathrm{Fr}^{\mathrm{crit}}$ is also constant.

The second term in the max-function in \eqref{eq:crcrit_orig} is not log-convex, because the exponent includes the decision variable $\mathrm{Fr}$.
We resolve this by introducing an auxiliary variable $r^\mathrm{Fr}$ and fitting a posynomial function to data generated from the power function:
    \begin{align} 
                 C_{R}^\mathrm{Frcrit} & \geq \max \left\{ 1.0, r^\mathrm{Fr} \right\}, \label{eq:crcrit1} \\
                 \left(r^\mathrm{Fr}\right)^{0.02608} & \geq 0.1528\mathrm{Fr}^{1.537}
    + 0.9672\mathrm{Fr}^{-0.008538}. \label{eq:crcrit2}
    \end{align}
The fit is excellent, with root mean square error less than 0.01\%.
\subsection{Cell and battery pack} \label{sec:batt_model}
We let $b_{ij}$ denote the battery discharging power on the sailing leg from port $i$ to $j$. The discharging power is the sum of shaft power and auxiliary power,
\begin{equation*}
    b_{ij} = P^\text{shaft}_{ij}+P^\text{aux}, \quad b_{ji} = P^\text{shaft}_{ji}+P^\text{aux},
\end{equation*}
for all $(i,j)\in \mathcal{L}$. Here $\mathcal{L}$ is the set of outbound sailing legs of service as described in §\ref{sec:indexing_sets}. The shaft power is the effective towing power divided by propeller open water efficiency: $P^\text{shaft}=(R_F+R_R)v/\eta^\text{prop}$.

The battery capacity must be sufficient to supply the energy demand of the sailing leg with the greatest energy consumption. The required energy supply is multiplied by the excess capacity factor $\Theta$ to impose a constraint on the minimum capacity $Q$ as
\begin{equation} \label{eq:battery_capacity}
    Q \geq \frac{\Theta}{\eta^\text{dis}} \max \left\{ E_{ij}^\text{leg} : (i,j)\in \mathcal{L}_s \right\},
\end{equation}
where $E_{ij}^\text{leg}=(bt^\text{transit})_{ij}/\eta^\text{dis}$ is the electrical energy drawn from the battery on the sailing leg $(i,j)$.

The battery is assembled from $N^\text{cell}$ number of individual lithium-ion cells with a specified fixed capacity $\tilde{Q}$ and constant voltage $V^\text{cell}$. The battery power $b$, cell current $\tilde{b}$, battery capacity $Q$ and cell capacity are related according to
\begin{equation*}
    \tilde{b} = \frac{b}{V^\text{cell}N^\text{cell}}, \quad \tilde{Q} = \frac{Q}{V^\text{cell}N^\text{cell}}.
\end{equation*}
We do not track the evolution of the cell charge over consecutive sailing legs. The mean cell charge normalized to its capacity is taken as 50\%, i.e., $\tilde{q}/\tilde{Q}=0.5$ with $\tilde{q}$ denoting the cell charge. This assumption always holds when the capacity factor $\Theta$ in \eqref{eq:battery_capacity} is equal to or greater than two.
\subsection{Battery capacity degradation}
As the battery cycles, its capacity decreases. The battery reaches the end of its lifetime once the capacity drops below some fraction of its initial value (e.g., 80\%). We adopt the semi-empirical cell degradation modeling approach in \cite{SO18}, which uses a combination of first-principles electrochemistry models and fitting of experimental data.

\begin{figure}
\begin{center}
\includegraphics[width=0.45\textwidth]{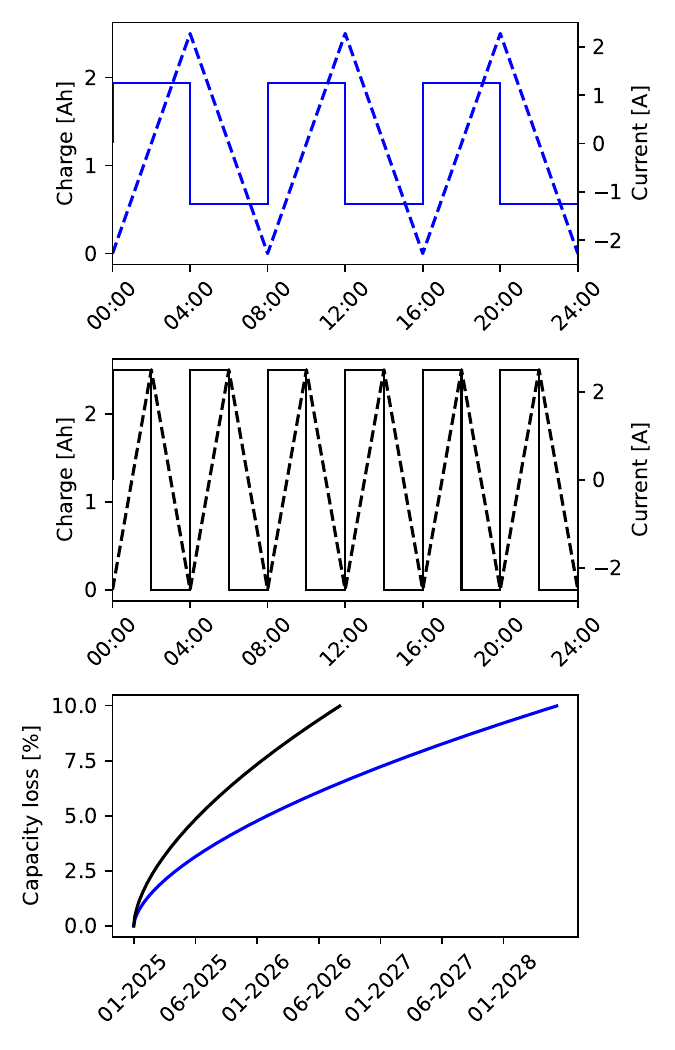}
\caption{Effect of cycle intensity on cell capacity loss. Three cycles per day (top), six cycles per day (middle), and capacity loss for both charging profiles (bottom).} 
\label{fig:cell_deg}
\end{center}
\end{figure}

The fractional capacity loss of a single cell (from its initial capacity) at time $t$ is
\begin{multline} \label{eq:cap_loss}
    \phi(t)= \\ \left[ \gamma^\text{acc}(t) \right]^{\chi_1} \left[ \chi_2 + \gamma^\text{cha}(t) \chi_3 \right] \exp \left( \frac{-E_A+\chi_4 \gamma^\text{cur}(t) }{R_GT^c} \right).
\end{multline}
The three terms in \eqref{eq:cap_loss} model the effects of accumulated charge, average charge, and average charge rate.
\begin{enumerate}
    \item In the first term, $\gamma^\text{acc}(t)=\int_0^t \left|\tilde{b}(\tau)\right| d\tau$ is the accumulated absolute charge throughput until time $t$.
    \item In the second term, $\gamma^\text{cha}(t)=(1/t)\int_0^t (\tilde{q}(\tau)/\tilde{Q}(\tau)) d\tau$ is the average normalized cell charge, which we take as 0.5. 
    \item The third term originates from the Arrhenius temperature dependence equation. Here $\gamma^\text{cur}(t)=(1/t)\int_0^t (\left|\tilde{b}(\tau)\right|/\hat{Q}(\tau)) d\tau$ is the average current rate until time $t$. The constants are activation energy $E_A$, universal gas constant $R_G$, and cell temperature $T^c$.
\end{enumerate}
The parameters $\chi_1,\ldots,\chi_4$ are identified from cell-dependent experimental data. We use the data from a 3.3V lithium-iron-phosphate cell in \cite{SO18}. Figure \ref{fig:cell_deg} illustrates the capacity loss of a single cell with varying cycle intensity.

Since we model the voyage as a set of sequential sailing legs with constant speed, the integrals in \eqref{eq:cap_loss} take the form of finite sums. We define the auxiliary term $r^\text{dis}$ and rewrite \eqref{eq:cap_loss} as the inequality 
\begin{equation} \label{eq:cap_loss_ineq}
        \phi^\text{max} \geq \left( N^\text{life}N^\text{rt}2 r^\text{dis} \right)^{\chi_1} \xi \exp \left( \frac{\chi_42N^\text{rt}r^\text{dis}}{R_GT^\text{c}\tilde{Q}t^\text{ph}} \right),
\end{equation}
where
\begin{equation*}
    \begin{aligned}
    \xi & = \left( \chi_2+\frac{1}{2}\chi_3 \right)\exp \left( \frac{-E_A}{R_GT^\text{c}} \right), \\ 
    r^\text{dis} &= \sum_{(i,j)\in \mathcal{L}_s} \tilde{b}_{ij}t_{ij}^\text{transit}.
    \end{aligned}
\end{equation*}
Here we set the maximum allowed capacity loss to a fixed value $\phi^\text{max}$ and define $N^\text{life}$ as a decision variable for the number of consecutive planning horizon periods that the battery is used in its lifetime. The constraint  \eqref{eq:cap_loss_ineq} is a valid log-convex inequality because the right-hand side is a product of a monomial and a log-convex exp-function.

We let $N^\text{batt}$ denote the number of complete batteries we need over the ship's lifetime. Since one battery is always installed for the newbuilding, $N^\text{batt}$ is defined as
\begin{equation} \label{eq:num_batt}
    N^\text{batt} \geq \max \biggl\{1, \frac{t^\text{life}}{N^\text{ph}N^\text{life}} \biggr\},
\end{equation}
where $N^\text{ph}$ is the number of planning horizon periods in a year and $t^\text{life}$ is the ship lifetime in years. 
The variable $N^\text{life}$ is used in the cost objective \eqref{eq:cost_objective} to model the battery investment capital cost. 
\subsection{Battery deck arrangement}
\begin{figure}
\begin{center}
\includegraphics[width=0.4\textwidth]{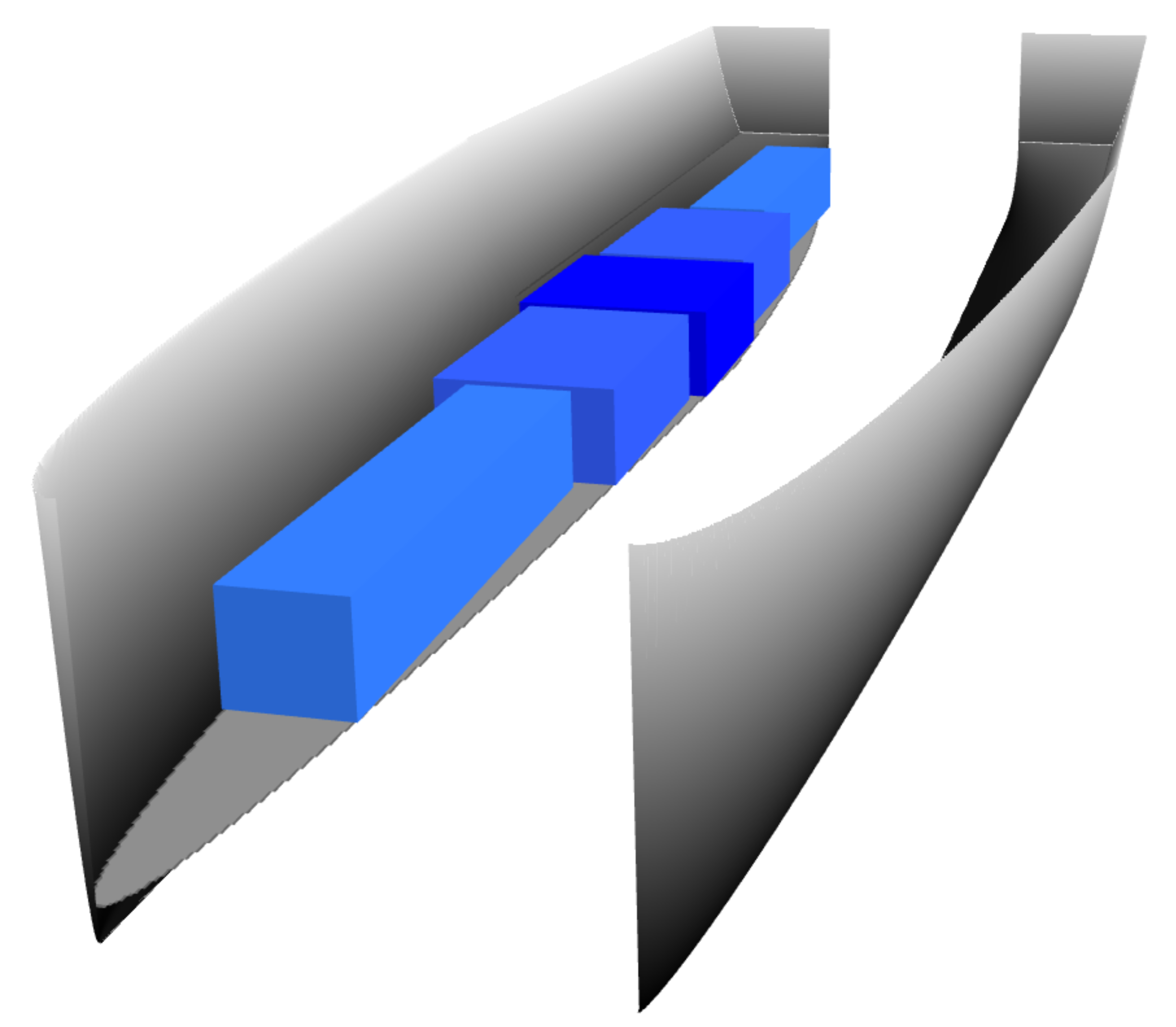}
\caption{The battery deck is divided into an uneven number of dedicated rooms whose locations and dimensions are design variables. The arrangement is symmetrical about the centermost room. The longitudinal center of gravity is constrained to lie above the hull center of buoyancy.} 
\label{fig:room_arrangement}
\end{center}
\end{figure}
The battery is assembled from many individual modules, stacked vertically in strings and housed in dedicated ventilated rooms. We consider $N^\text{room}$ uneven number of rooms with the same capacity $Q^\text{room}=Q/N^\text{room}$. The arrangement is symmetrical about the centermost room (figure \ref{fig:room_arrangement}).

The \textit{k}th room is specified by its width $w_k^\text{batt}$, length $l_k^\text{batt}$, equal height $h^\text{batt}$, and lower corner $(w_k/2, \tilde{l}_k, \tilde{h} )^\text{batt}$ in relation to the reference point located at the keel centerline of the forward perpendicular. We impose lower bounds on the room height and elevation from the keel:
\begin{align} 
    h^\text{batt} & \geq h^\text{batt,min}, \label{eq:room_lower_bound1} \\ 
    \tilde{h}^\text{batt} & \geq 0.1D. \label{eq:room_lower_bound2}
\end{align}

The room dimensions relate to the volume via monomial equations
\begin{equation*}
    (w_khl_k)^\text{batt} = Q^\text{room}\rho^\text{vol}
\end{equation*}
for each $k=1,\ldots,\lceil N^\text{room}/2\rceil$ with $k=1$ denoting the foremost room. The parameter $\rho^\text{vol}$ is the volumetric energy density of a battery space.

Each room is enclosed within the hull and prevented from overlapping with adjacent rooms:
\begin{align} 
        &\frac{w_k^\text{batt}}{2} \leq H^\text{fore}(\tilde{l}_k^\text{batt}, \tilde{h}^\text{batt}),\, k = 1,\ldots,\lceil N^\text{room}/2\rceil, \label{eq:room_positioning1} \\
        &l_k^\text{batt} + \tilde{l}_k^\text{batt} \leq \tilde{l}_{k+1}^\text{batt},\quad k = 1,\ldots,\lfloor N^\text{room}/2\rfloor. \label{eq:room_positioning2}
\end{align}
A ship is designed to float on an even keel by longitudinally aligning the ship's weight and the buoyant force vectors. The position of the centermost room determines the battery center of gravity because the room arrangement is symmetrical. The resulting longitudinal stability constraint is
\begin{equation} \label{eq:long_force_balance}
   \tilde{l}^\text{batt}_{k=\lceil N^\text{room}/2\rceil} + \frac{1}{2} l^\text{batt}_{k=\lceil N^\text{room}/2\rceil}  \leq \tilde{l}^\text{cb}.
\end{equation}
The inequality is tight because it forces the room position forward from its volume-maximizing position.  

The buoyant force acts at the centroid of the displaced water volume. The longitudinal distance of the centroid from the forward perpendicular, denoted $\tilde{l}^\text{cb}$ in \eqref{eq:long_force_balance}, is given by the first moment integral
    \begin{equation*}
        \begin{aligned}
                \tilde{l}^\text{cb} & = \frac{1}{\nabla} \int_{0}^{L/2} \int_{0}^{T} (yH^\text{fore}(y,z) + (y+1)H^\text{aft}(y,z))\ dz\ dy  \\ 
                & = \left[ \frac{2\beta^3/(5\hat{\beta}) + 3\beta^2/2+ \beta - \tilde{\beta}\log{\hat{\beta}}}{2\beta^3/(3\hat{\beta}) + \beta^2 - \beta \log \hat{\beta}} \right]L,
        \end{aligned}
    \end{equation*}
where the parameters $\hat{\beta}$ and $\tilde{\beta}$ are the same as in \eqref{eq:centroid_vert_pos}.

The relative positioning constraints of the decks are
\begin{align} 
    h^\text{batt} + \tilde{h}^\text{batt} & \leq \tilde{h}^\text{roro}, \label{eq:roro_spacing1} \\  
    T & \leq \tilde{h}^\text{roro}, \label{eq:roro_spacing2} \\ 
    \tilde{h}^\text{roro} + h^\text{roro} & \leq D. \label{eq:roro_spacing3}
\end{align}
The first two constraints state that the first ro-ro deck is above the battery deck and the waterline. The third constraint specifies ship depth as the distance from keel to the second ro-ro deck, which we take as the uppermost continuous deck.
\subsection{Weight breakdown} \label{sec:weight_breakdown}
Lightweight and deadweight are calculated as
\begin{align}
        W_L &= \underbrace{\rho^\text{st}\sum_k N_k^\text{pl}A_k^\text{pl}p_k}_\text{Structures} + \underbrace{(0.3L + \left( L^\text{sup} \right)^{3.1}/10^5)B }_\text{Outfitting} \nonumber \\
         & + \underbrace{\rho^\text{mass}Q/\rho^\text{vol}}_\text{Battery} +  \underbrace{\rho^\text{mot}P^\text{shaft, max}}_\text{Motors}, \label{eq:weight_breakdown1} \\
        W_D &= \underbrace{1565q^\text{cap}_{c=\text{roro}} }_\text{Ro-ro cargo} + \underbrace{170(0.02+1.146q^\text{cap}_{c=\text{pax}})}_\text{Passengers and crew} \nonumber \\
        & + \underbrace{170(0.02+1.146q^\text{cap}_{c=\text{pax}})}_\text{Freshwater}. \label{eq:weight_breakdown2}
\end{align}
Outfitting weight calculation uses adjusted coefficients from \cite{WG77}. The basis of calculation for steel weight is the actual or estimated area and thickness of each steel plate element. Table \ref{tb:plate_elements} presents posynomial formulas for area calculations and the quantity and thickness of each element. The notation $A_k^\text{tbh}$ refers to the area of the \textit{k}th transverse bulkhead (extending from keel to battery deck ceiling) defined as
\begin{equation*}
    \begin{aligned}
      A_k^\text{tbh} &= \int_0^{\tilde{h}^\text{roro}} H^\text{fore}(\tilde{l}_k^\text{batt},z)\,dz \\ 
      &=  \frac{\beta}{\beta+1}\tilde{h}^\text{roro}_kB \left(\frac{\tilde{h}^\text{roro}_k}{T} \right)^{1/{\beta}}\sqrt{\frac{2\tilde{l}^\text{batt}_k}{L}}.
    \end{aligned}
\end{equation*}
\begin{table}
\small
\begin{center}
\caption{Surface areas and thicknesses of plate elements.}\label{tb:plate_elements}
\begin{tabular}{llll} \hline
Plate element & $N^\text{pl}$ & $A^\text{pl}$ & $p^\text{pl}$ \\ \hline
External hull & 2 & $\varphi_ALl^\mathrm{mid}_D$ & $p^\text{sp}$ \\ 
Bottom plate & 1 & $\frac{1}{2}LB$ & $p^\text{bp}$ \\
Inner bottom & 1 & $\frac{2}{3}\left(\frac{\tilde{h}^\text{batt}}{T} \right)^{1/\beta} BL $ & $p^\text{td}$ \\ 
Ro-ro decks & 2 & $\frac{5}{6}BL$ & $p^\text{td}$ \\ 
Ro-ro ramp & 1 & $Bh^\text{roro}$ & $p^\text{sup}$ \\
Longitudinal& 2 & $0.7LD$ & $p_\text{side}$ \\
bulkheads  &&& \\
Transverse & 2 & $\sum_{k=1}^{\lceil N^\text{room}/2\rceil}A_k^\text{tbh}$ & $p^\text{sup}$ \\ 
bulkheads &&& \\
Walls & 2 & $\Bigl( 2\sum_{k=1}^{\lfloor N^\text{room}/2\rfloor}l_k^\text{batt} $ & $p^\text{sup}$ \\ 
& & $ + \, l_{k=\lceil N^\text{room}/2\rceil}^\text{batt}  \Bigr)h^\text{batt} $ & \\
Superstructure & 1 & $L^\text{sup}(6h^\text{sup} + 4B + 2h^\text{roro})$ & $p_\text{sup}$ \\ 
& & $+ \, B(h^\text{roro} + 6h^\text{sup}) $ &  \\ \hline
\end{tabular}
\end{center}
\end{table}
\subsection{Cost terms}
The total annualized ship acquisition and variable operation cost (excluding energy) in \eqref{eq:obj} is
    \begin{equation} \label{eq:cost_objective}
        \begin{aligned}
            C^\text{ship} & = C^\text{hotel}f^\text{cap}_\text{pax} + C^\text{deck} +  C^\text{batt} N^\text{batt} Q +  C^\text{hull} V^\text{GT}.
        \end{aligned}
    \end{equation}
The cost terms are hotel crew, deck crew, battery, and construction.

The construction cost and port charges are levied in proportion to the ship's gross tonnage $V^\text{GT}$ calculated as 
\begin{equation} \label{eq:gross_tonnage}
    \begin{aligned}
        V^\text{GT} & = V^\text{int}(0.2 + 0.02\log_{10}(V^\text{int})) \approx V^\text{int} 0.2 \\
        & + 0.02\log_{10}(\hat{V}^\text{int}) V^\text{int} \left( \frac{V^\text{int}}{\hat{V}^\text{int}} \right)^{1/\log_{10}(\hat{V}^\text{int})\log(10)},
    \end{aligned}
\end{equation}
where $V^\text{int}$ is the internal volume of the ship:
\begin{equation*}
    V^\text{int} = C^\text{int}_BLBD + N^\text{sup}h^\text{sup}L^\text{sup}.
\end{equation*}
(Here $C^\text{int}_B$ is an internal volume block coefficient.) The last term in \eqref{eq:gross_tonnage} replaces the log-function with its local monomial approximation at $\hat{V}^\text{int}$. A reasonable choice of $\hat{V}^\text{int}$ is $10^5$. The error from the approximation is less than 0.1\% for typical ship sizes.
\section{Computational experiments} \label{sec:comp_experiments}
This section focuses on benchmarking the nonlinear formulation of the ESFSMP (with continuous ship sizing) against the conventional linear formulation (with fixed ship types). The aim is to assess the practical value of developing and employing the posynomial models that enable using a physics-based continuous ship sizing in ESFSMP. We formulate test cases of varying size and solve both nonlinear and linear formulations using the same open-source and commercial solvers for convex mixed-integer nonlinear programming problems. The nonlinear and linear problem formulations are compared in terms of the resulting problem sizes, i.e, the number of decision variables and constraints. Computational performance is compared in terms of total cost objective value and solution time. 
\subsection{Input data}
\begin{figure*}
\begin{center}
\includegraphics[width=0.9\textwidth]{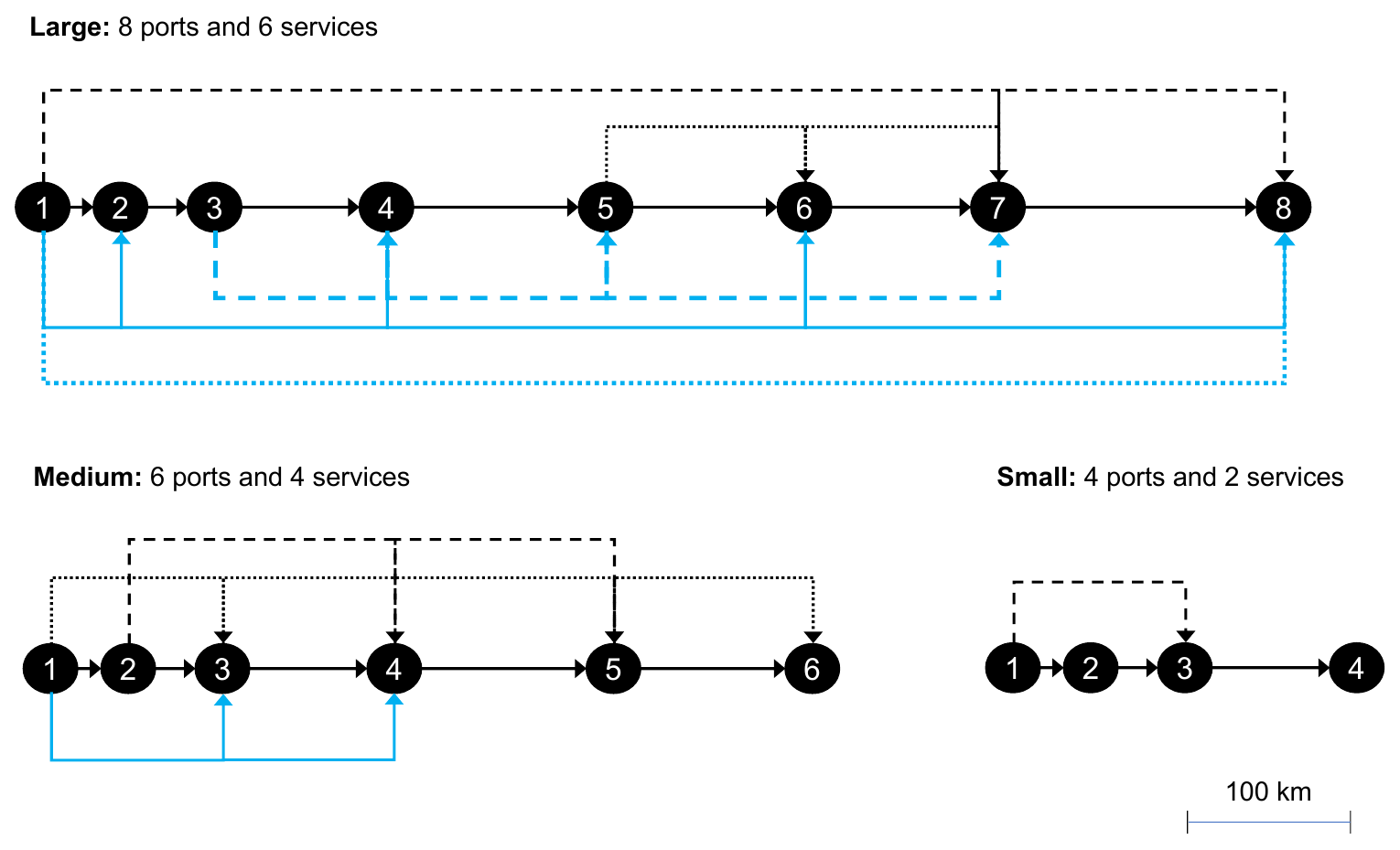}
\caption{Randomly generated shipping networks. The directed arcs of different styles denote the outbound sailing legs of each service.} 
\label{fig:ex_networks}
\end{center}
\end{figure*}

We construct a shipping network consisting of eight ports ordered on a straight line. The distances between the ports are randomly generated from a uniform distribution on the interval $[40, 200]$ (km). The demands of ro-ro and passenger traffic are also randomly generated from a uniform distribution on the interval $[0.5, 5]$ (1000 passengers or 1000 lane-meters). A scenario labeled \textit{small} consists of two pendulum services that use the first four ports. Two other scenarios labeled \textit{medium} and \textit{large} consist of four services over the first six ports and six services over all the ports, respectively. The first service in each scenario visits every port. The other services are randomly generated. Figure \ref{fig:ex_networks} illustrates the ports and the outbound sailing legs of each service. 
\subsection{Linear benchmark problem}
\ref{sec:lin_prob} formulates the linear ESFSMP. The modeling principle, based on \cite{HLN+22}, is to define a small number of discrete speed and energy levels and linearly interpolate the other values. Each discrete level is associated with a weight variable. The linear model also introduces binary variables for selecting ship types from a predetermined set of alternatives. 

There are two additional minor restrictions compared to the ESFSMP formulated in \ref{sec:logcvx_formulation}. First, every port is equipped with the same fixed charger power. Second, cargo offloading and loading time at each port are independent of the cargo volume. We apply these additional restrictions to the nonlinear problem instances as well to ensure that the performance comparisons are valid. 
\subsection{Results}
We formulate linear and nonlinear problem instances for the three input data scenarios described above. The implementations are available via the link provided in §\ref{sec:introduction}. The linear instances use 21 ship types, which are generated using the ship model from §\ref{sec:vessel_model}. Table \ref{tb:prob_size} lists the number of integer variables, continuous variables, and constraints in each problem instance. The modeling software transforms the nonlinear (log-convex) instances to a standard power cone form, which adds a large number of continuous auxiliary variables. The linear instances contain at least 20x the number of integer variables compared to the nonlinear instances. 

We solve all the problem instances using MOSEK 11.0, which is a high-performance commercial solver for convex and mixed-integer optimization problems. We also use ECOS 2.0 \citep{DCB13} to solve the nonlinear instances and SCIP 10.0 \citep{HBB+25} to solve the linear instances. Both are free to use non-commercial solvers. The computing platform is a single-core 2.25 GHz processor with 13.3 GB of memory. 

Table \ref{tb:lin_solves} reports the solution times and the optimal objective function values for the linear problem instances. MOSEK solves the small instances almost 3x faster than SCIP, but is more than 1.6x slower solving the medium instance and 2.4x slower solving the large instance. This reflects the typical performance variability of mixed-integer optimizers. The main message from table \ref{tb:lin_solves} is that the problem quickly becomes intractable when the number of integer variables increases. The worst-case solution time is less than one and a half minutes for the small instance but over three hours for the large instance.

Table \ref{tb:nonlin_solves} reports the results for the nonlinear problem instances. Total cost objective value improves by 11.7 to 19.7\%, while the worst-case solve time is reduced by 10.5 to 19.6x compared to the linear instances. The non-commercial solver ECOS is capable of solving all the instances and shows worse performance than MOSEK only for the large problem instance.  

When interpreting the results, it is important to recognize that the solvers find the global optimum (within a small tolerance) of all the instances. The differences in solution quality between the linear and nonlinear instances do not originate from poor local solutions. Instead, the differences are attributed to the fact that the nonlinear formulation can find improved ship designs outside of the 21 predetermined types in the linear formulation.

\begin{table*}
\small
\begin{center}
\caption{Problem size as reported by MOSEK after the pre-solve stage.}\label{tb:prob_size}
\begin{tabular}{llllll} \\\hline
Problem  & Formulation  & Integer variables  & Continuous variables & Constraints  \\  \hline
Small      & Linear     & 477  & 830  & 1079   \\ 
Small      & Nonlinear  & 20    & 1544  & 612     \\ 
Medium     & Linear     & 881   & 1775 & 2691     \\ 
Medium     & Nonlinear  & 40    & 4675 & 1831    \\ 
Large      & Linear     & 1410  & 3328  & 5023      \\ 
Large      & Nonlinear  & 56    & 7900 & 3058     \\ 
\hline
\end{tabular}
\end{center}
\end{table*}
\begin{table}
\small
\begin{center}
\caption{Solver performance for the linear problem instances.}\label{tb:lin_solves}
\begin{tabular}{llll} \\\hline
Problem    & Solver & Objective     & Solving \\ 
           &        & [MEUR]        & time [s] \\ \hline
Small      & MOSEK  & 68.2         & 22.6 \\ 
Small      & SCIP   & 68.2         & 81.4 \\ 
Medium     & MOSEK  & 98.2         & 2181.3 \\ 
Medium     & SCIP   & 98.2         & 1317.2 \\ 
Large      & MOSEK  & 227.2         & 10804.3 \\ 
Large      & SCIP   & 227.2        & 4495.7 \\ 
\hline
\end{tabular}
\end{center}
\end{table}
\begin{table}
\small
\begin{center}
\caption{Solver performance for the nonlinear problem instances.}\label{tb:nonlin_solves}
\begin{tabular}{llll} \\\hline
Problem    & Solver & Objective & Solving \\ 
           &        & [MEUR]    & time [s] \\ \hline
Small      & ECOS   & 54.8     & 1.6 \\ 
Small      & MOSEK  & 54.8     & 7.7 \\ 
Medium     & ECOS   & 86.8     & 44.4 \\ 
Medium     & MOSEK  & 86.8     & 129.3 \\ 
Large      & ECOS   & 198.3    & 551.2 \\ 
Large      & MOSEK  & 198.3    & 223.5 \\ 
\hline
\end{tabular}
\end{center}
\end{table}
\section{Baltic Sea ferry service design}  \label{sec:example_problem}
\subsection{Transportation network}
We model sea transportation services in the northern Baltic Sea region. The route network connects five major ports in Sweden, Estonia, the Finnish mainland, and the autonomous region of the Åland Islands. Connections include both short-distance shuttle services and overnight cruises. Most ships are ro-pax ferries offering transportation for both ro-ro cargo and passengers. The total number of passengers transported in all the routes was 11.8 million in 2024 \citep{Tal25}. Ro-ro cargo volume, consisting of passenger vehicles, trucks, and trailers, was 2.2 million in the same period.

Conventional ships currently operate all the routes in this region. Fast bunkering of liquid hydrocarbon fuels enables short turnaround times. In particular, the 160 nautical mile route from Turku (Finland) to Stockholm (Sweden) (figure \ref{fig:routes}, top left) runs on a 24-hour periodic schedule with less than an hour turnaround time in both end ports. Battery-electrification of this route using the current schedule requires charging powers exceeding 100 megawatts or implementing a battery-swapping scheme.

The other routes' relaxed schedules or short distances make them better suited for electrification. The long-haul overnight routes from Helsinki and Tallinn to Stockholm run on a 48-hour periodic schedule with 7-hour port visits. The flexibility in these routes creates an opportunity to optimize sailing speeds or add a new port call to the service. For example, the Helsinki–Mariehamn–Stockholm ships can add a 5-hour shuttle service to Tallinn instead of staying berthed continuously for 7 hours at Helsinki.
\subsection{Cases}
We formulate and solve six instances of the ESFSMP in §\ref{sec:logcvx_formulation}, which are labeled as
\begin{equation*}
    \text{B4, U4, M4, M3, M2 and M1}.
\end{equation*}
The letter on the label stands for the fleet type. In the baseline fleet B, the number of ships and their principal dimensions match the present fleet of the largest ferry operator in the region. Moreover, we assign the services the same frequencies currently used in the routes. The uniform fleet U consists of interchangeable ships of the same optimized design allocated for all the routes. The mixed fleet M is assembled from ships optimized for each route. 

\begin{figure*}
\begin{center}
\includegraphics[width=1.0\textwidth]{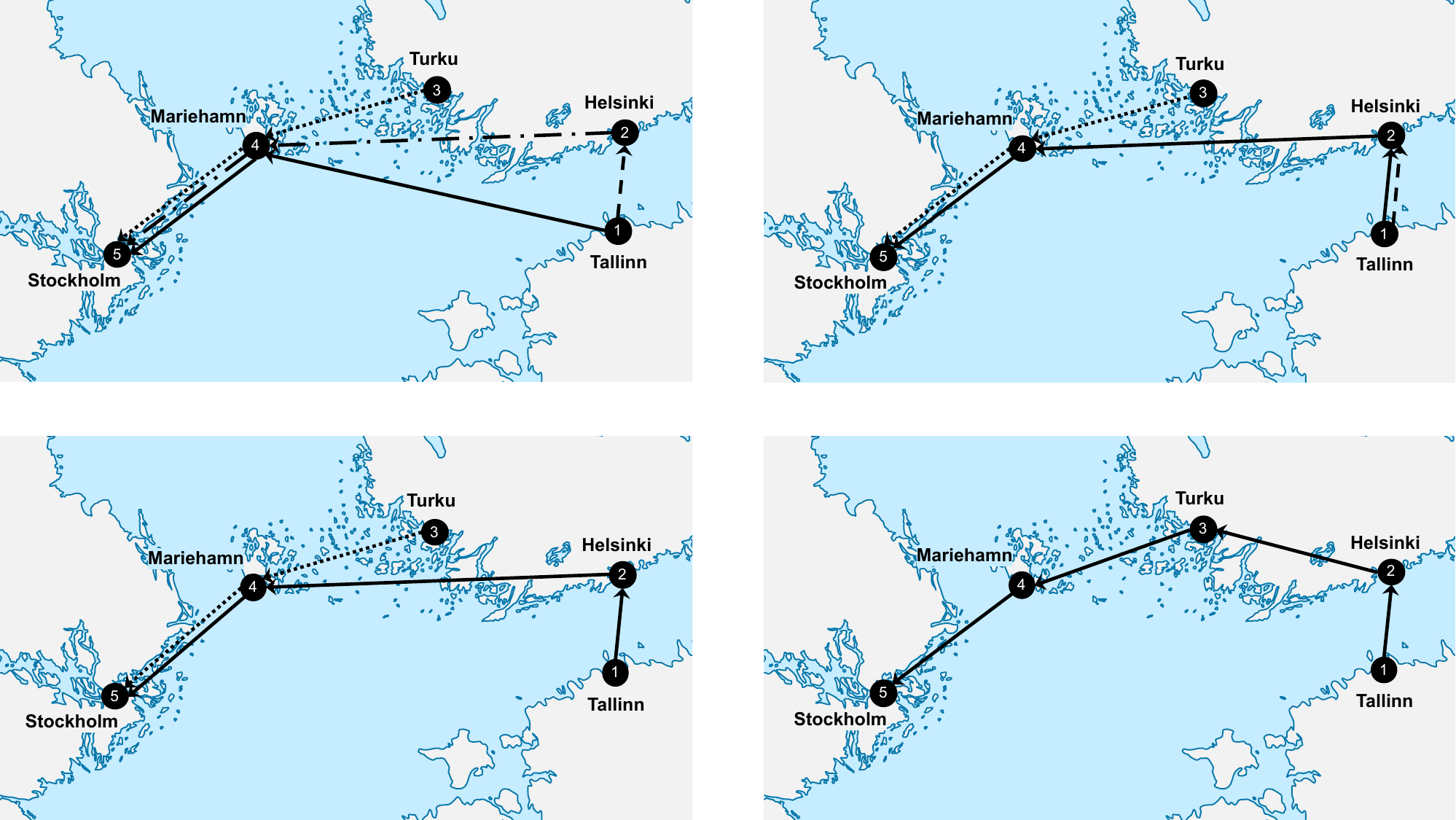}
\caption{Four alternative route plans for a northern Baltic Sea ferry service. A route is a sequence of sailing legs denoted by lines of the same style. A complete voyage consists of outbound and inbound sailing along the same
set of legs.} 
\label{fig:routes}
\end{center}
\end{figure*}

The number in the label represents the number of distinct routes in the transportation network. Figure \ref{fig:routes} (top left) illustrates the current network of four routes. The first alternative network in figure \ref{fig:routes} (top right) routes traffic from Estonia to Sweden via Finland, eliminating one route. The number of routes is further reduced to two and one while preserving a direct connection between all the mainland ports (bottom panels in figure \ref{fig:routes}). In all the alternative networks, ships in service to and from Sweden call at the Åland Islands to benefit from tax exemption, which enables duty-free sales on board.

Table \ref{tb:demands} shows the demand input data (the same for all problem instances) for every active origin-destination port pair. The other fixed parameter values can be accessed in the implementation code. We solve all the problem instances using ECOS with the default settings.
\begin{table}
\small
\begin{center}
\caption{48-hour demands for each active port-pair. The data are given in the format \textit{passengers (cargo)}. One unit is a thousand lane-meters or a thousand passengers. The demand is estimated from route-specific traffic volumes and market shares reported in \cite{Tal25}}\label{tb:demands}
\begin{tabular}{llllll} \hline
$i/j$ 	& 1 		& 2 		& 3 		& 4        & 5  \\ \hline
1 		& - 		&  20.50 	& - 		& 0.05     & 1.72 \\
  		& - 		&  (47.92) 	& - 		& (0.01)   & (3.00)  \\
2 		& 20.50 	&  - 		& - 		& 0.1      & 4.56  \\
  		& (47.92) 	&  - 		& - 		& (0.02)   & (2.21) \\
3 		& - 		&  - 		& -	 		& 1.36     & 6.10 \\
  		& - 		&  - 		& - 		& (0.25)   & (19.85) \\
4 		& 0.05  	&  0.1  	& 1.36 	    & -	       & 1.71  \\
  		& (0.01) 	& (0.02) 	& (0.25) 	& - 	   & (0.32) \\
5 		& 1.72 	    & 4.56 	    & 6.10 	    & 1.71     & - \\
  		& (3.00) 	& (2.71) 	& (19.85) 	& (0.32)   & - \\ \hline
\end{tabular}
\end{center}
\end{table}
\subsection{Results}

\begin{table*}
\small
\begin{center}
\caption{Optimal values of selected fleet design variables. The annual total fleet cost in millions of EUR is displayed in brackets below the case label. Sailing leg speeds are in the order they appear in the outbound voyage. The value in the brackets is the speed of the same sailing leg in the inbound voyage.}\label{tb:solution}
\begin{tabular}{lllllllllllllll} \hline
Case    & Route & $N^\text{trip}$      & $N^\text{ship}$        & $L$          & $L^\text{sup}$ & $V^\text{GT}$     & $E^\text{batt}$ & $t^\text{cell}$  & $v_1$ & $v_2$ & $v_3$ & $v_4$ \\ 
        &       &            &              & [m]          & [m]            & $10^3$   & [MWh] & [a]      & [kn] & [kn] & [kn] & [kn] \\ \hline
B4      & 12    & \textbf{18} & \textbf{3}   & \textbf{212} & \textbf{191}   & 43.6     & 38.9 & 2.1       & 21.6 (21.6) & -  & -  & - \\
(216.4) & 145   & \textbf{2} & \textbf{2}   & \textbf{212} & \textbf{191}   & 55.3     & 82.8 & 7.6       & 14.9 (15.0) & 15.6 (15.3) & -  & -  \\
        & 245   & \textbf{2} & \textbf{2}   & \textbf{203} & \textbf{183}   & 50.4     & 81.6 & 7.6       & 15.1 (15.1) & 15.6 (15.5) & - & -  \\
        & 345   & \textbf{2} & \textbf{1}   & \textbf{212} & \textbf{191}   & 58.9     & 62.3 & 2.5       & 15.7 (15.6) & 15.7 (15.8) & - & - \\ \hline
U4      & 12    & 6          & 1            & 161          & 130            & 21.9     & 101.6 & 17.0     & 16.0 (16.0) & - & - & - \\
(114.5) & 145   & 2          & 1            & 161          & 130            & 21.9     & 101.6 & 4.2      & 21.2 (21.3) & 23.1 (23.7) & - & - \\
        & 245   & 4          & 2            & 161          & 130            & 21.9     & 101.6 & 4.6      & 21.2 (21.2) & 23.1 (22.8) & - & -  \\
        & 345   & 6          & 3            & 161          & 130            & 21.9     & 101.6 & 16.5     & 14.3 (14.2) & 14.3 (14.2) & - & - \\ \hline
M4      & 12    & 6          & 1            & 109          & 83             & 11.1     & 13.1 & 1.6      & 15.1 (15.1) & - & -  & - \\
(106.6) & 145   & 2          & 1            & 200          & 101            & 30.3     & 131.6 & 4.3     & 21.8 (21.7) & 22.7 (22.8) & - & - \\
        & 245   & 2          & 1            & 223          & 202            & 54.4     & 163.6 & 4.5     & 22.3 (22.3) & 23.3 (23.2) & -  & - \\
        & 345   & 6          & 3            & 133          & 119            & 18.7     & 24.8 & 3.2      & 14.2 (14.1) & 14.4 (14.4) & - & - \\ \hline
M3      & 12    & 4          & 1            & 96           & 63             & 9.9      & 5.7 & 3.6       & 10.0 (10.0) & - & - & -  \\
(103.8) & 1245  & 4          & 2            & 185          & 167            & 32.8     & 197.1 & 4.2     & 27.6 (27.2) & 25.6 (25.6) & 27.2 (27.0) & - \\
        & 345   & 4          & 2            & 141          & 127            & 21.8     & 24.5 & 3.9      & 13.9 (13.9) & 14.6 (14.6) & - & -  \\ \hline
M2      & 1245  & 4          & 2            & 210          & 150            & 36.8     & 224.2 & 4.3     & 27.7 (27.3) & 25.9 (25.9) & 27.3 (27.1) & - \\ 
(120.7) & 345   & 6          & 3            & 143          & 120            & 22.0     & 24.8 & 3.9      & 13.9 (13.8) & 14.5 (14.5) & - & - \\ \hline
M1      & 12345 & 14          & 14           & 166          & 107            & 20.5     & 60.3 & 6.7      & 17.5 (17.5) & 16.1 (16.1) & 17.5 (17.5) & 17.5\\ 
(193.8) & & & & & & & & & & & & (17.5) \\  \hline
\end{tabular}
\end{center}
\end{table*}
\begin{table}
\small
\begin{center}
\caption{Optimal shore charger powers in megawatts.}\label{tb:charge_powers}
\begin{tabular}{llllll} \hline
Case    & 1 & 2 & 3 & 4 & 5 \\ \hline
B4      & 8.9 & 8.9 & 21.2 & 77.5 & 21.0 \\
U4      & 75.5 & 59.1 & 8.6 & 131.0 & 56.7 \\
M4      & 79.2 & 89.9 & 9.47 & 164.3 & 54.0 \\
M3      & 27.1 & 133.6 & 11.2 & 182.8 & 49.3  \\
M2      & 23.5 & 149.1 & 9.5 & 202.2 & 42.7 \\
M1      & 3.0 & 19.1 & 29.8 & 26.7 & 5.2 \\ \hline
\end{tabular}
\end{center}
\end{table}
Table \ref{tb:solution} shows optimal values of selected fleet design and operation variables, while table \ref{tb:charge_powers} shows the optimal port charger powers. Comparing the solutions for cases B4 and M4, an optimized fleet more than halves the fleet cost compared to the present fleet in the same 4-route network. M4 deploys smaller ships to the short routes and a smaller overall fleet size. On the other hand, two routes operate at a higher frequency, enabled by higher speeds and shorter turnaround times. In the Helsinki-Tallinn route, M4 maintains the same frequency as B4, but manages to reduce speed by optimizing the port turnaround times. The ship deployed to the cargo-heavy Tallinn-Helsinki route features a short superstructure, which saves building cost, electricity, and port charges.  

Deploying a uniform fleet to the 4-route network is only 7.4\% more expensive than an optimized fleet (cases U4 and M4). This result is surprising considering that the gross tonnage of the largest ship in M4 is almost five times larger than the gross tonnage of the smallest ship. In U4, the batteries are vastly oversized in the short routes. However, the oversizing also brings low cycling load and long lifetime.

The 3-route network modeled in the case M3 achieves the lowest cost among all the cases, which indicates that the currently used 4-route network is suboptimal. The basic idea in the 3-route network is to reduce fleet size by combining the routes originating from Tallinn and Helsinki to Stockholm. Simultaneously, the new combined route brings additional capacity to the Tallinn-Helsinki route. The port of Tallinn is equipped with a smaller charger because the frequency of the Helsinki-Tallinn drops from three daily round-trips to two. Finally, the single and 2-route networks M1 and M2 show poor performance because they supply reduced capacity to the high-demand Helsinki-Tallinn route.
\subsection{Parameter sensitivity}
\begin{figure}
\begin{center}
\includegraphics[width=0.43\textwidth]{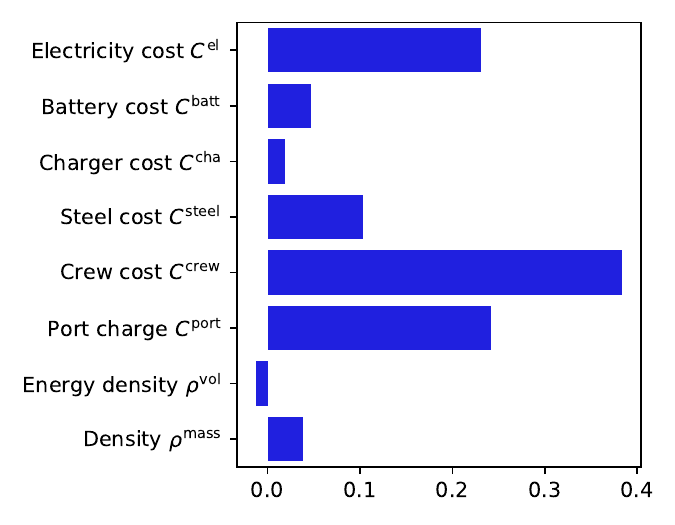}
\caption{Sensitivity of selected parameters for the case M3.} 
\label{fig:param_sens}
\end{center}
\end{figure}

Figure \ref{fig:param_sens} reports the tractional sensitivity of the objective function to fractional changes to selected fixed parameters for the case M3. A sensitivity of +1 means that perturbing the parameter value by 1\% increases the objective value by 1\% in the vicinity of the optimal solution. The sensitivities are practically useful for guiding the planning process to focus on the most important strategic decisions. 

All the parameter sensitivities are nonzero. Therefore, at least one constraint that includes a given parameter is active at the optimum, governing the solution. Interestingly, the two most sensitive parameters (crew cost and port charge) are unrelated to the battery-electric powertrain. The sensitivity of battery density is positive, because battery mass per unit of energy increases linearly with density. 
\section{Extensions and variations}
\subsection{Demand scenarios}
Demand for maritime transportation services varies daily, weekly, and monthly. For example, the monthly passenger traffic volume in the northern Baltic Sea is three times greater in the summer holiday season compared to mid-winter \citep{statfin25}. Incorporating demand scenarios into the problem formulation in §\ref{sec:logcvx_formulation} is straightforward by introducing indexing of the operational decision variables by the demand scenarios. This extension costs many new decision variables and constraints, including the discrete service frequency variables. However, the problem decomposes into smaller scenario-specific subproblems, only coupled by the fleet design variables. Specialized solution methods can leverage this property.  
\subsection{Flow aggregation}
The log-sum function \eqref{eq:constr5} assigns a different utility for a transportation volume depending on the number of services that deliver that volume between an origin-destination port pair. Specifically, the formulation favors a larger number of services. An alternative approach is to calculate the utility from the aggregate flow from all the services. The constraints that aggregate the flows impose a lower bound on posynomials, which destroys the favorable log-convex structure. In \cite{BKV+07} (§9.1), the authors suggest a sequential local solution method that involves taking local monomial approximations of the lower-bounded posynomials. They observe that this method works well in practice.
\subsection{Travel time cost}
Thus far, we have ignored that the demand between a given origin-destination port pair is not fixed, but depends on the transportation service's availability and characteristics. Passenger ferry services compete with low-cost airlines, and ro-ro cargo can use (at least in some cases) alternative land-based routes. To ensure that the new battery-electric ferry fleet offers competitive travel times, we can extend the cost-objective \eqref{eq:cost_objective} with a travel time term. Since the total travel time on a route is a posynomial expression of sailing leg speeds, the travel time terms preserve the favorable log-convex structure.
\section{Conclusions}
The design of a battery-electric liner shipping service presents a complex optimization problem due to the strong coupling of decisions concerning cargo routing, sailing schedules, charger powers, and the ship sizing subdomains. Many ship owners will need to deploy battery-electric fleets in the near future, as increasingly strict emission regulations are enforced. The ro-ro, ro-pax, and container feeder segments constitute over 1,350 ships in total in European Union waters alone \citep{Ep25}. As the transition to an electric marine transportation system involves a clean-sheet design of all aspects of the system, there is a unique opportunity to apply rigorous mathematical optimization techniques to make a large impact.

This paper has presented a convex formulation of the battery-electric liner shipping fleet size and mix problem. The problem is to minimize total fleet purchase and operation cost while modeling everything from fleet size down to the hull girder plate thickness of each ship. A typical real-world numerical case is represented by a network consisting of six ports, four services, and 30 origin-destination demand arcs. The solve time is less than two minutes on a desktop computer, demonstrating the novel possibility for fast fleet design optimization by employing computationally efficient convex optimization methods.
\section*{Acknowledgement}
This work was supported by the Future Shipping Electrified (FUSE) research project through Business Finland under grant 3545/31/2023.
\appendix
\section{Changes of variables} \label{sec:var_changes}
A nonlinear optimization problem defined by log-convex posynomial functions is written as
    \begin{equation} \label{eq:gp_prob}
        \begin{aligned}
            & \text{minimize} && \sum_{k=1}^{K_0}c_k\prod_{l=1}^{n}x_l^{a_{lk}} \\
            & \text{subject to} && \sum_{k=1}^{K_m}c_k\prod_{l=1}^{n}x_l^{a_{lk}} \leq 1, \quad m = 1,...,M. \\
        \end{aligned}
    \end{equation}
In the standard form \eqref{eq:gp_prob}, the objective is a nonconvex function and the constraints define nonconvex sets. However, the problem has a special structure that can be exploited. A logarithmic change of variables $y_j=\log x_j$ (such that $x_l = \exp y_l$) allows a transformation that convexifies \eqref{eq:gp_prob}. Specifically, the log of a posynomial is convex in $y$:
\begin{equation} \label{eq:log-exp-sum}
    \begin{aligned}
    &\log\sum_{k=1}^{K}c_k\prod_{l=1}^{n}x_l^{a_{lk}} = \log\sum_{k=1}^{K}c_k\exp{(a_k^Ty)} \\
    &= \log\sum_{k=1}^{K}\exp{(\log c_k + a_k^Ty)}
    \end{aligned}
\end{equation}
because log-sum-exp is a convex function \citep{BV04} and the argument of the exponential function is linear in $y$. 

Using the log transformed posynomials of the form \eqref{eq:log-exp-sum}, and the abbreviated constant $b_k = \log c_k$, the log-convex problem is reformulated as the convex optimization problem
    \begin{equation} \label{eq:gp-cvx}
        \begin{aligned}
            & \text{minimize} && \log\sum_{k=1}^{K_0}\exp{(b_{0k} + a_{0k}^Ty)} \\
            & \text{subject to} && \log\sum_{k=1}^{K_m}\exp{(b_{mk} + a_{mk}^Ty)} \leq 0, \quad m = 1,...,M. \\
        \end{aligned}
    \end{equation}
The optimal solution of the original nonconvex problem \eqref{eq:gp_prob} is easily recovered from a solution of the transformed problem \eqref{eq:log-exp-sum} since $x_l^* = \exp y_l^*$.
\begin{figure}
\begin{center}
\includegraphics[width=0.31\textwidth]{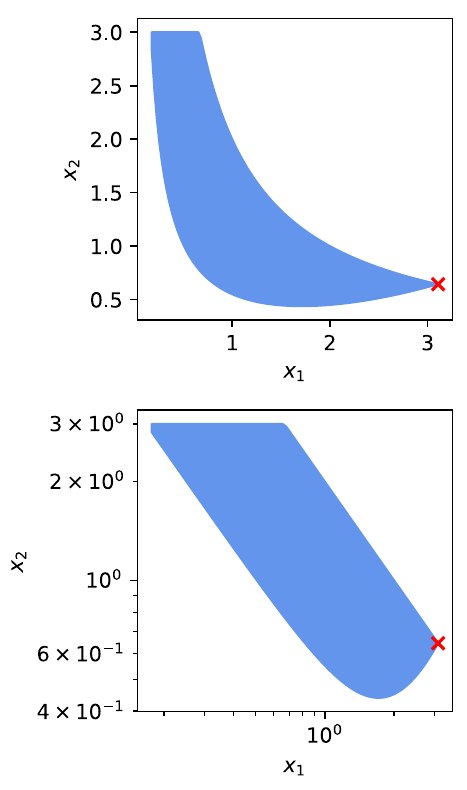}
\caption{Feasible set of a log-convex optimization problem with two variables. The top panel shows the original nonconvex problem, and the bottom panel shows the log-transformed convex problem. The solid red cross is the optimal solution.} 
\label{fig:constr_set}
\end{center}
\end{figure}

The following example illustrates convexification via log-transformation. The problem is to minimize a monomial subject to a lower bounding posynomial and two upper bounding monomials with variables $x = (x_1, x_2) \in \mathbb{R}_{++}^2$:
        \begin{equation*}
        \begin{aligned}
            &\text{minimize} \quad && x_1^{-1} \\
            &\text{subject to} \quad && x_2 \geq 0.5x_1^{-1} + x_1^2/20 \\
            & && x_2 \leq 2x_1^{-1},\, x_2 \leq 3.
        \end{aligned}
    \end{equation*}
Figure \ref{fig:constr_set} draws the feasibility set of the problem in both linear and log-space.
The red cross is the global optimum $x^* \approx (3.11, 0.64)$.
The constraints define a nonconvex set in linear space but a convex set in log-space.
\section{Still water loads} \label{sec:sw_loads}
\begin{figure}
\begin{center}
\includegraphics[width=0.4\textwidth]{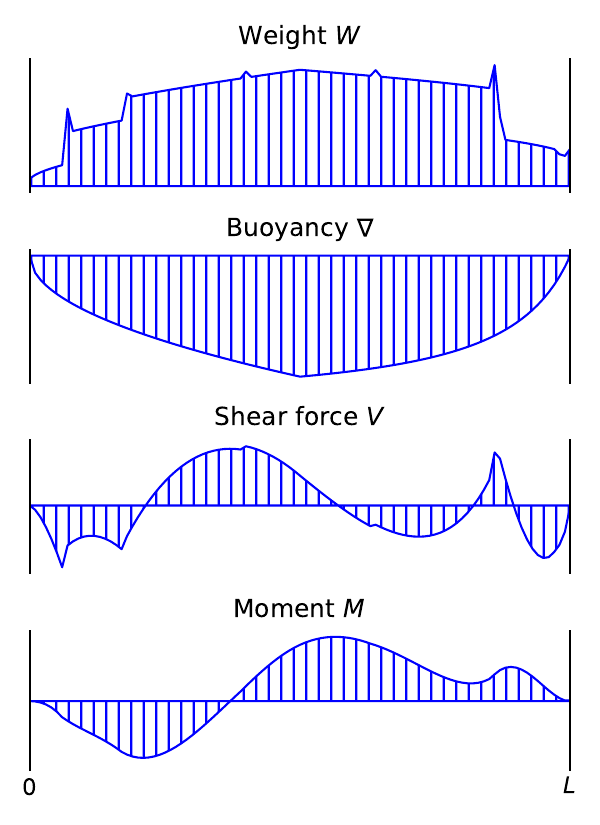}
\caption{Weight, buoyancy forces, and the resulting still water loads. The origin of the coordinate system is at the forward perpendicular.} 
\label{fig:sw_loads}
\end{center}
\end{figure}

The top two panels in figure \ref{fig:sw_loads} illustrate buoyancy and weight per unit length of the hull in §\ref{sec:weight_buo_balance} floating with static draft at even keel. The local buoyancy force per unit length equals the local hull submerged cross-sectional area. The weight distribution is the sum of distributions of the weight components \eqref{eq:weight_breakdown1} and \eqref{eq:weight_breakdown2}, assuming $L^\text{sup}/L=0.8,$ $N^\text{room}=3$ and $l^\text{batt}_k/L=0.23$. The middle panel in figure \ref{fig:sw_loads} shows the imbalance of the distributions along the length of the hull. This imbalance causes shear force and bending moment at each hull cross-section. 

Let $W(y)$ and $\nabla(y)$ denote the unit length weight and buoyancy at the longitudinal coordinate $y$. The still water shear force $V^\text{sw}(y)$ and bending moment $M^\text{sw}(y)$ lengthwise distributions are given by the integrals
\begin{equation*}
    V^\text{sw}(y) = \int_0^y (W(\tilde{y})-\nabla(\tilde{y}))d\tilde{y}, \quad M^\text{sw}(y) = \int_0^y V^\text{sw}(\tilde{y})d\tilde{y}.
\end{equation*}
The resulting forces and moments per unit length are illustrated in the bottom two panels in figure \ref{fig:sw_loads}.

\section{Hull girder section modulus}  \label{sec:sec_modulus}
Normal stress $\sigma(z)$ at a vertical distance $z$ from the reference line relates to the bending moment $M$ as
    \begin{equation*}
            \sigma(z) = \frac{M}{\tilde{I}}(z-\tilde{z}_\text{NA}),
    \end{equation*}
where $\tilde{I}$ is a geometrical property of the cross-section, called the \textit{moment of inertia} about the girder neutral axis $\tilde{z}_\text{NA}$.
The section modulus $Z=\tilde{I}/(z-\tilde{z}_\text{NA})$ in relation to the structural element with the highest stress defines the ultimate strength of the hull girder.

The section modulus is a function of the principal dimensions and thicknesses of the structural elements of the hull girder illustrated in figure \ref{fig:cross_section}. The structural components are the bottom, inner bottom, first ro-ro deck, topmost continuous deck (i.e., second ro-ro deck), side shell, and longitudinal bulkheads. We model the two adjacent longitudinal bulkheads as a single bulkhead with double thickness.

The thicknesses of the side shell plates, lumped bulkhead, and bottom relate to the topmost deck according to the fixed ratios 
\begin{equation*}
    p^\text{sp} = \frac{Bp^\text{td}}{2D}, \quad p^\text{bh}=2p^\text{sp}\quad \text{and} \quad p^\text{bp}=\frac{2}{3} p^\text{td}.
\end{equation*}
The vertical position of the cross-section neutral axis from the reference is the vertical coordinate of the cross-section centroid, defined as
\begin{equation*}
    \tilde{z}_\text{NA}=\frac{\sum_iN_i(z_\text{NA})_iA_i}{\sum_iN_iA_i} = \frac{2}{5}D,
\end{equation*}
where the structural elements are indexed by $i$. The quantity $N_i$, neutral axis coordinate $(z_\text{NA})_i$ and area $A_i$ of each element are given in table \ref{tb:normal_stress_elements}.  

The second moment of area of the cross-section is calculated using the parallel axis theorem according to
\begin{equation} \label{eq:cross_section_2nd_moment}
    \begin{aligned}
        \tilde{I}=&\sum_i \left( I_i+A_i(z_\text{NA}-\tilde{z}_\text{NA})_i^2 \right) \\
        =&\sum_i \left( \frac{w_ih_i^3}{12}+w_ih_i(z_\text{NA}-\tilde{z}_\text{NA})_i^2 \right) \\
        \approx & \frac{133}{150}p^\text{td}BD^2,
    \end{aligned}
\end{equation}
where $w_i$ is the width of the \textit{i}th element, $h_i$ is the height, and $I_i=w_ih_i^3/12$ is the second moment of area of the element with respect to its own centroid. In \eqref{eq:cross_section_2nd_moment}, we use the approximation $I_i\approx0$ for all the deck and bottom elements. The second moment of area of a thin-walled horizontal element, about its own neutral axis, is sufficiently small to be negligible.
\begin{table}
\small
\begin{center}
\caption{Properties of the midship structural elements.}\label{tb:normal_stress_elements}
\begin{tabular}{lllllll} \hline
Element     & $N$   & $w$               & $h$            & $A$             & $z_\text{NA}$   & $z_\text{NA}$  \\ 
 & &  &     &         &   & $-\tilde{z}_\text{NA}$ \\ \hline
Bottom      & 1     & $B$               & $\frac{3p^\text{td}}{2}$ & $\frac{3Bp^\text{td}}{2}$ & 0               & $\frac{-2D}{5}$  \\
Inner bottom    & 1     & $B$               & $p^\text{td}$            & $Bp^\text{td}$            & $\frac{D}{10}$ & $\frac{-3D}{10}$ \\ 
Ro-ro deck      & 1     & $B$               & $p^\text{td}$            & $Bp^\text{td}$            & $\frac{D}{2}$  & $\frac{D}{10}$ \\ 
Top deck       & 1     & $B$               & $p^\text{td}$            & $Bp^\text{td}$            & $D$             & $\frac{3D}{5}$ \\ 
Side plate & 2     & $p^\text{sp}$   & $D$            & $\frac{Bp^\text{td}}{2}$  & $\frac{D}{2}$  & $\frac{D}{10}$ \\ 
Bulkhead    & 1     & $p^\text{sp}$    & $D$            & $Bp^\text{td}$            & $\frac{D}{2}$  & $\frac{D}{10}$ \\ \hline
\end{tabular}
\end{center}
\end{table}

Since $\tilde{z}_\text{NA}<D/2$, the distance from the cross-section neutral axis to the girder top deck is larger than the distance to the bottom.
The section modulus at the top deck is lower than at the bottom, which implies that the maximum normal stress is at the deck.
The deck section modulus is
    \begin{equation*}
            Z^\text{deck} = \frac{5\tilde{I}}{3D}.
    \end{equation*}
\section{Hull girder shear flow} \label{sec:shear_flow}
The vertical shear force $V$ induces a shear stress $\tau$ in the hull girder. The shear force is equal to the shear stress integrated over the area of the cross-section. The shear stress distribution in a thin-walled cross-section is determined using \textit{shear stress flow} $q=p\tau$. We let $s$ denote a coordinate running along the cross-section. Shear flow per unit shear force at $s$ is
\begin{equation} \label{eq:shear_flow_integration}
    q(s)=q_0-\frac{1}{\tilde{I}}\int_0^s(z(s)-\tilde{z}_\text{NA})p(s)\,ds,
\end{equation}
where $q_0$ is the shear flow at the starting point of the integration \citep{VRR17}. 

\begin{figure}
\begin{center}
\includegraphics[width=0.4\textwidth]{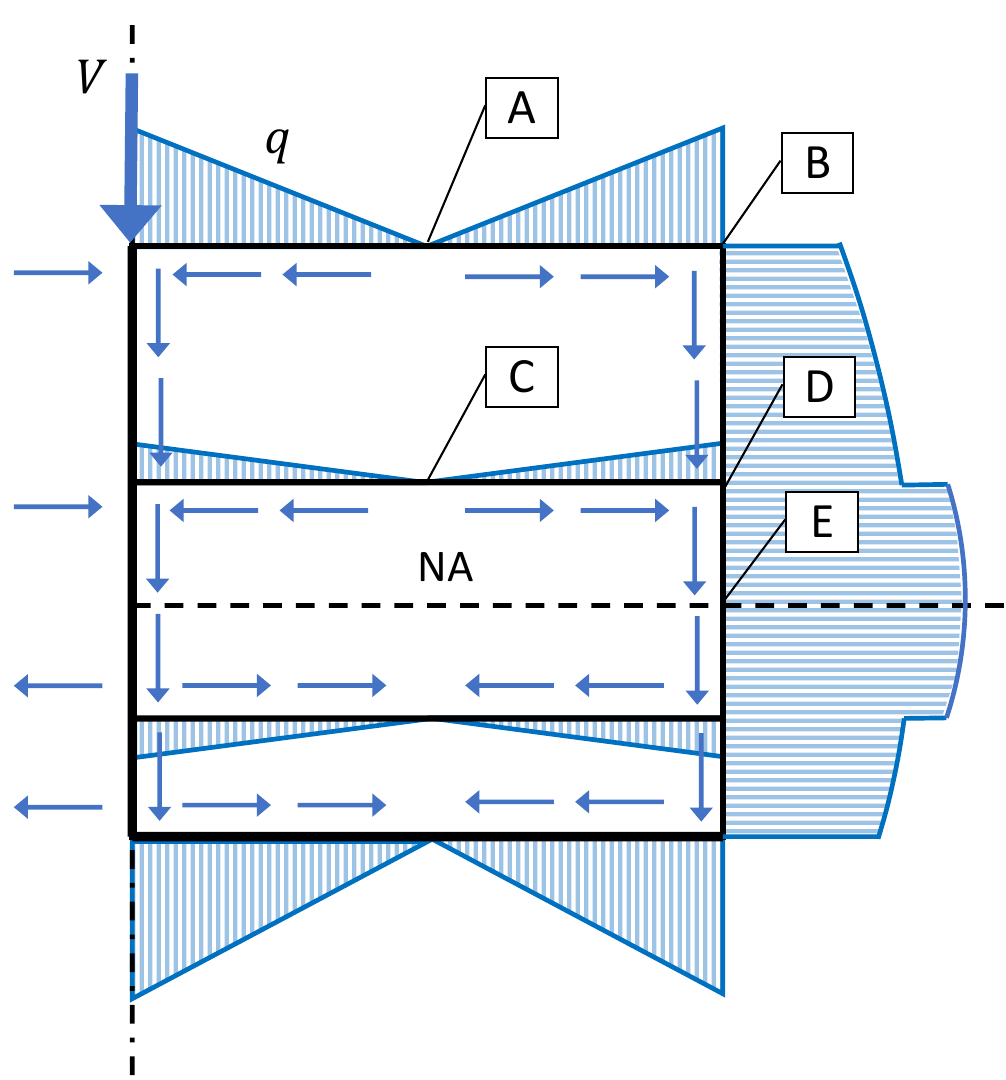}
\caption{Hull girder half-section shear flow distribution. The maximum shear flow value is located at the neutral axis. The value at the bulkhead is twice that of the side shell plates. Since the bulkhead is also twice as thick as the side shell plates, the maximum shear stress values are equal.} 
\label{fig:shear_flow_dist}
\end{center}
\end{figure}

The coordinates where $q(s)=0$ are indeterminate in a cross-section of one or more closed cells. When the starting point of the shear flow integration in \eqref{eq:shear_flow_integration} (with $q_0=0$) is arbitrary, discontinuities in the longitudinal displacement $u(s)$ may occur. The discontinuities need to be compensated by adding constant shear flows along the walls of all the closed cells until the integration of the longitudinal displacement is zero along the circumference of each cell:
\begin{equation*}
   \oint_\text{cell} u(s)\,ds = - \frac{1}{G} \oint_\text{cell} \frac{q(s)}{p(s)} \,ds=0.
\end{equation*}
(Here $G$ is the shear modulus.) For the details of this compensation procedure, see §5 in \cite{VRR17} or §6.3 in \cite{Mit21}. The symmetry properties of the hull girder cross-section result in zero shear flows at the half-section vertical centerline of each horizontal element, as illustrated in figure \ref{fig:shear_flow_dist}.

The hull girder cross-section is symmetrical with respect to the skeleton line of the bulkhead, and the bulkhead is twice as thick as the side shell plates. Therefore, the same maximum shear stress value occurs at the neutral axis of both elements. In the following, we use the maximum shear flow value at the side shell plates to determine the ultimate stress.

The shear flow \eqref{eq:shear_flow_integration} at the end node $j$ of a constant thickness straight line segment from node $i$ to $j$ is 
\begin{equation} \label{eq:shear_flow_line_seg}
    q_j = -\frac{pl}{2\tilde{I}}(z_j+z_i-2\tilde{z}_\text{NA})+q_i,
\end{equation}
where $z_i,z_j$ are the $z$-coordinates (from baseline) of the nodes $i,j$, the length of the line segment is $l$ and $q_i$ is the shear flow at node $i$ \citep{lloyds23}. Using nodes from A to E (figure \ref{fig:shear_flow_dist}) and four line segments defined in table \ref{tb:line_segments}, we find the maximum shear flow per unit shear force
\begin{equation*}
    q_E = -\frac{53BDp^\text{td}}{200\tilde{I}}
\end{equation*}
as a function of the midship cross-section principal dimensions and the deck plate thickness. 
\begin{table}
\small
\begin{center}
\caption{Properties of the hull girder cross-section line segments.}\label{tb:line_segments}
\begin{tabular}{llllllll} \\\hline
Element     & $i$ & $j$ & $l$ & $p$ & $z_i$ & $z_j$ & $q_i$ \\ \hline
Top deck   & A & B & $\frac{B}{4}$ & $p^\text{td}$ & D & D & 0 \\
Side shell  & B & D & $\frac{D}{2}$ & $p^\text{sp}$& D & $\frac{D}{2}$ & $q_{A \rightarrow B}$ \\
Ro-ro deck   & C & D & $\frac{B}{4}$ & $p^\text{td}$ & $\frac{D}{2}$ & $\frac{D}{2}$ & 0 \\ 
Side shell   & D & E & $\frac{D}{10}$ & $p^\text{sp}$ & $\frac{D}{2}$ & $\frac{2D}{5}$ & $q_{A \rightarrow D}$ \\ 
&&&&&&& $+q_{C \rightarrow D}$ \\ \hline
\end{tabular}
\end{center}
\end{table}
\section{Estimation of wetted area} \label{sec:wetted_area}
The wetted area is the surface integral
    \begin{equation*}
            A_S = 2 \int_{0}^{L/2} \int_{0}^{T} \left( \left\Vert \begin{bmatrix}
                -\frac{\partial H^\text{fore}}{\partial y} \\ 
                -\frac{\partial H^\text{fore}}{\partial z} \\
                1
            \end{bmatrix}
            \right\Vert + 
            \left\Vert \begin{bmatrix}
                -\frac{\partial H^\text{aft}}{\partial y} \\ 
                -\frac{\partial H^\text{aft}}{\partial z} \\
                1
            \end{bmatrix}
            \right\Vert
            \right) \ dz\ dy,
    \end{equation*}
where $\left\Vert \cdot \right\Vert$ denotes the 2-norm.
However, we cannot evaluate the above integral precisely.
Instead, we use the approximation
\begin{equation} \label{eq:wet_area_app}
      A_S \approx 2\varphi_ALl^\mathrm{mid}_T,
\end{equation}
where $l^\mathrm{mid}_T$ is the midship arc length from keel centerline to hull waterline and $\varphi_A$ is a correction factor that accounts for the narrowing hull towards the forward and aft of amidships. The arc length is
\begin{equation} \label{eq:mid_arc_length}
    \begin{aligned}
    l^\mathrm{mid}_T &= \int_{0}^{B/2} \left\Vert \begin{bmatrix}
                -\frac{d z^\text{mid}(x)}{d x} \\
                1
            \end{bmatrix}
            \right\Vert\, dx \\
            &= \int_{0}^{B/2} \sqrt{1+\left(\frac{2}{B} \right)^{2\beta}(T\beta)^2x^{2\beta-2} } \, dx,
    \end{aligned}
\end{equation}
where $z^\text{mid}(x)$ is z-coordinate of the hull at amidships at offset $x$. By writing $x=H^\text{fore}(L/2,z)$ from \eqref{eq:hull_form1} with respect to $z$, we obtain
$z^\text{mid}(x) = T\left(2x/B\right)^\beta$.

The integral \eqref{eq:mid_arc_length} can be evaluated exactly for only a few $\beta$ values. We deal with the general case by substituting the integrand function with an interpolating polynomial. A cubic polynomial is sufficiently accurate because the integrand is smooth and curves only in one direction. Integration using a cubic polynomial (known as Simpson's 3/8 rule) yields
\begin{equation} \label{eq:simp_app}
    l^\mathrm{mid}_T \approx \frac{B}{16}\sum_{k'=0}^{3}S_{k'+1}g^\mathrm{mid}(k'B/6),
\end{equation}
where $g^\text{mid}$ is the integrand function in \eqref{eq:mid_arc_length} and $S=(1,3,3,1)$ is the vector of Simpson's multipliers.

The integrand $g^\mathrm{mid}$ evaluated in \eqref{eq:simp_app} is a posynomial raised to a fractional power and thus not log-convex. Next, we introduce auxiliary variables and an inequality constraint to obtain an equivalent log-convex form. Using the auxiliary variables $r^g_{k'},\; k'=0,\ldots,3$ we can now write \eqref{eq:wet_area_app} as the pair of valid posynomial inequalities
\begin{align}
        A_S & \geq  2L\varphi_A\sum_{k'=0}^{3}S_{k'+1}\sqrt{r^g_{k'}},  \label{eq:wet_area1} \\
        r^g_{k'} & \geq 1+\left(\frac{2}{B} \right)^{2\beta}(T\beta)^2 \left( \frac{k'B}{6}\right)^{2\beta-2},  \quad k'=0,\ldots,3. \label{eq:wet_area2}
\end{align}
\section{Linear problem formulation} \label{sec:lin_prob}
We need two new sets for the mathematical description of the linear problem: let $\mathcal{U}$ be the set of ship types indexed by $u$ and $\mathcal{V}_{u}$ be the set of discrete speed levels for ship type $u$ indexed by $v$. 

The energy consumption of a ship type $u$ sailing at speed $v$ from port $i$ to $j$ is a fixed parameter $E_{ijvu}$. Ship type $u$ has capacity $q^\text{cap}_{cu}$ of cargo $c$. The ship purchase cost $C^\text{ship}_{u}$ and port charge $C^\text{port}_{iu}$ are also ship type dependent. The transit time from port $i$ to $j$ depends on the speed level $v$ and is represented by $t^\text{transit}_{ijv}$. The cargo handling time $t^\text{cargo}_{ij}$ is fixed at the destination port $j$. We use binary variables to select one frequency per service from available frequencies represented by $N^\text{trip}_{f}$.


We define the following new decision variables: $x_{ijsvu}\in[0,1]$ is a weight for speed level $v$ for ship type $u$ on service $s$ on leg $(i,j)$. As the weights sum to one on each leg, service, and ship type, the true speed is a linear combination of the discrete speed levels. The charging time at the destination port of leg $(i,j)$ in service $s$ for ship type $u$ is $t^\text{cha}_{ijsu}\in\mathbb{R}_{++}$ and the total roundtrip duration in service $s$ is $t^\text{trip}_{s}\in\mathbb{R}_{++}$. The ship types and service frequencies are selected using binary variables. Let $z_{fs}\in\{0,1\}$ equal 1 if frequency $f$ is chosen on service $s$ and 0 otherwise and similarly let $y_{us}\in\{0,1\}$ equal 1 if ship type $u$ is selected and 0 otherwise. The nonnegative discrete variable $N^\text{ship}_{su}\in \mathbb{Z}_{\geq 1}$ is the number of ships of type $u$ deployed on service $s$. In addition, we need auxiliary variables $q_{ijcs}^\text{logapp}\in\mathbb{R}$ for the piecewise linear approximation of the log utility function of $q_{ijcs}$.

The problem formulation is as follows:
\phantom{Some text here that flows to the next line a bit more}
minimize
\begin{equation}
    \begin{aligned}
        & \sum_{s \in \mathcal{S}} \sum_{u \in \mathcal{U}} \Biggl[C^\text{ship}_uN^\text{ship}_{su}
        + \sum_{(i, j) \in \mathcal{L}_s} \sum_{f \in \mathcal{F}} N^\text{trip}_fz_{sfu} C_{ju}^\text{port} \\
        &+ \sum_{(i, j) \in \mathcal{L}_s} \sum_{f \in \mathcal{F}} C_j^\text{el}N^\text{trip}_fz_{sfu}t^\text{cha}_{ijs}P_j \Biggr],
         \label{eq:lin_obj}
    \end{aligned}
\end{equation}
subject to
\begin{align}
& N^\text{ship}_{su} \leq N^\text{ship,max}_uy_{us},\quad s\in \mathcal{S}, \, u\in \mathcal{U}, \label{eq:lin_constr1} \\
& \sum_{u \in \mathcal{U}} y_{us} = 1, \quad s\in \mathcal{S}, \label{eq:lin_constr2} \\
& \sum_{f \in \mathcal{F}} z_{sfu} =  y_{us},\quad s\in \mathcal{S}, \, u\in \mathcal{U} ,\label{eq:lin_constr3} \\
& 
\begin{array}{lr}
    t^\text{trip}_{s} = \sum_{(i,j) \in \mathcal{L}_s} \Bigl( \sum_{v \in \mathcal{V}_u} \sum_{u \in \mathcal{U}} t^\text{transit}_{ijv}x_{ijsvu} \\ + t^\text{port}_{ijs}\Bigr), \quad
    s\in\mathcal{S},
\end{array}
\label{eq:lin_constr4} \\
& \sum_{u \in \mathcal{U}} \sum_{f \in \mathcal{F}} N^\text{trip}_fz_{sfu} t_{s}^\text{trip} \leq \sum_{u \in \mathcal{U}} N^\text{ship}_{su}t_s^\text{route}, \quad s\in \mathcal{S}, \label{eq:lin_constr5} \\
& t^\text{cha}_{ijs} = \frac{1}{P_j} \sum_{v \in \mathcal{V}_u}\sum_{u \in \mathcal{U}}E_{ijvu}x_{ijsvu}, \, (i,j) \in \mathcal{L}_s,\, s\in \mathcal{S},\label{eq:lin_constr6} \\
& t_{ijs}^\text{port} \geq \max \{ t_{ij}^\text{cargo}, t_{ijs}^\text{cha} \}, \quad (i,j) \in \mathcal{L}_s,\, s\in \mathcal{S},\label{eq:lin_constr7} \\
& \sum_{v \in \mathcal{V}_u} x_{ijsvu} = y_{us}, \quad (i,j) \in \mathcal{L}_s,\, s\in \mathcal{S}, \, u\in \mathcal{U}, \label{eq:lin_constr8} \\
& \sum_{v \in \mathcal{V}_u}E_{ijvu}x_{ijsvu} \leq E^\text{max}_uy_{us},
\begin{array}{ll}
& (i,j) \in \mathcal{L}_s,\\ 
& s\in \mathcal{S}, \, u\in \mathcal{U},
\end{array} \label{eq:lin_constr9} 
\\
& 
\begin{array}{ll}
    \sum\limits_{i'\in \mathcal{N}^{-}_{is}} \sum\limits_{j'\in \mathcal{N}^{+}_{is} \backslash \{i\}} q_{i'j'cs} \leq  \sum\limits_{u \in \mathcal{U}}\sum\limits_{f \in \mathcal{F}} q_{cu}^\text{cap}N^\text{trip}_fz_{sfu}, \\
    c\in \mathcal{C}, \, s\in \mathcal{S},\, i\in \mathcal{N}_s \backslash \{d_s\},
\end{array} \label{eq:lin_constr10} \\
& \begin{array}{ll}
    \sum\limits_{i'\in \mathcal{N}^{-}_{is}} \sum\limits_{j'\in \mathcal{N}^{+}_{is} \backslash \{i\}} q_{j'i'cs} \leq \sum\limits_{u \in \mathcal{U}}\sum\limits_{f \in \mathcal{F}} q_{cu}^\text{cap}N^\text{trip}_fz_{sfu}, \\ 
    c\in \mathcal{C}, \, s\in \mathcal{S},\, i\in \mathcal{N}_s \backslash \{d_s\},
\end{array}  \label{eq:lin_constr11} \\
& \sum_{s:(i,j) \in \mathcal{A}_s} q_{ijcs} \leq q_{ijc}^\text{dem},\quad c\in \mathcal{C},\ (i,j) \in \mathcal{A}, \label{eq:lin_constr12} \\
& \sum_{ (i,j) \in \mathcal{A} } \sum_{s:(i,j) \in \mathcal{A}_s} \sum_{c\in \mathcal{C}} \alpha_{ijc} q^\text{logapp}_{ijcs} \geq U^\text{min},\label{eq:lin_constr13} \\
& \begin{array}{lr}
         q^\text{logapp}_{ijcs} \leq 2\log(\hat{n})q_{ijcs}+\log(n) + 2n\log(1/\hat{n}), \\
         \quad (i,j) \in \mathcal{L}_s,\,c\in \mathcal{C}, \, s\in \mathcal{S}, \, n=1/2, 1, 2, 4,\ldots
    \end{array} \label{eq:lin_constr14}
\end{align}
where $\hat{n}=(n+1/2)/n$.

The objective function \eqref{eq:lin_obj} minimizes the total annualized cost. The three terms are the ship purchase and running cost, port charges, and electricity cost. In contrast to §\ref{sec:logcvx_formulation}, the shore charger installation cost is excluded because charger powers are fixed. 

The constraints are as follows:
\begin{enumerate}
    \item Constraints \eqref{eq:lin_constr1}-\eqref{eq:lin_constr3} encode logic relations: a nonzero number of ships of a given type can be deployed on a service only if the ship type is selected. Each service uses only one ship type and one frequency. 
    \item Constraints \eqref{eq:lin_constr4}-\eqref{eq:lin_constr8} couple service frequency and time spent on transit, port cargo handling, and charging. The total round-trip time \eqref{eq:lin_constr4} is the sum of speed-dependent sailing times and port turnaround times over all the legs. Ships must complete all their round-trips within the available time of the planning horizon as enforced by \eqref{eq:lin_constr5}. The charging time at the destination port of a leg is a function of the energy consumption on the leg by equation \eqref{eq:lin_constr6}. The port operation that takes the longest determines the turnaround time by \eqref{eq:lin_constr7}. The speed weights of the selected ship type sum to one in the equation \eqref{eq:lin_constr8}.
    \item The inequality \eqref{eq:lin_constr9} states that the energy consumption on any leg cannot exceed the battery capacity.
    \item The capacity constraints \eqref{eq:lin_constr10}-\eqref{eq:lin_constr11} ensure that the sum of cargo quantity that flows along a sailing leg does not exceed ship capacity. The constraints for the outbound and inbound parts of the voyage are given separately. The demand constraints \eqref{eq:lin_constr12} state that the total cargo quantity transported from an origin port to a destination port by all the services does not exceed the demand. The inequality formulation implies that demand can be rejected.
    \item The transportation service must provide a minimum level of utility by \eqref{eq:lin_constr13}. The constraints \eqref{eq:lin_constr14} encode a piecewise linear approximation of the nonlinear log utility function. 
\end{enumerate}

Lastly, we reformulate the expression $z_{sfu}t^\text{cha}_{ijs}$ in the objective \eqref{eq:lin_obj} and the expression $z_{sfu} t_{s}^\text{trip}$ in the constraints \eqref{eq:lin_constr4}. These terms consist of products of a binary variable and a nonnegative continuous variable. The idea is to replace the bilinear product term with a new continuous variable that is forced to take a positive value only when the binary variable is one. 

To reformulate the product term in \eqref{eq:lin_obj}, we introduce new auxiliary variables $\tilde{t}_{sfu}^\text{cha}\in\mathbb{R}_+$ (for all $s\in\mathcal{S},\,f\in\mathcal{F},\,u\in\mathcal{U}$), and the following new constraints:
\begin{equation} \label{eq:lin_reform}
   \tilde{t}_{sfu}^\text{cha}  \leq \sum_{(i,j)\in\mathcal{L}_s} t_{ijs}^\text{cha}-M^{(1)}(1-z_{sfu}), \quad s\in\mathcal{S},\,f\in\mathcal{F},\,u\in\mathcal{U}.
\end{equation}
Here $M^{(1)}$ is a sufficiently large constant. The last term in \eqref{eq:lin_obj} is then rewritten linearly as
\begin{equation*}
    \sum_{f \in \mathcal{F}} C_j^\text{el}N^\text{trip}_f\tilde{t}_{sfu}^\text{cha}P_j.
\end{equation*}
The constraint \eqref{eq:lin_reform} states that the left-hand side must be equal to or greater than the sum on the right-hand side if the binary variable is one. Otherwise, the left-hand side is reduced to the lower bound of zero.

The product terms in \eqref{eq:lin_constr4} are reformulated using the same technique as above.

\end{document}